\documentclass[12pt]{article}

\usepackage{amsfonts}
\usepackage{graphicx}
\usepackage{amsmath}
\usepackage{amsthm}

\newtheorem{theorem}{Theorem}[section]

\newtheorem{conjecture}[theorem]{Conjecture}
\newtheorem{corollary}[theorem]{Corollary}

\newtheorem{lemma}[theorem]{Lemma}

\newtheorem{question}[theorem]{Question}
\newtheorem{proposition}[theorem]{Proposition}

\newtheorem{fact}[theorem]{Fact}
\newtheorem{protoeg}[theorem]{Example}
\newtheorem{protoremark}[theorem]{Remark}
\newtheorem{protodefinition}[theorem]{Definition}

\renewcommand{\thetheorem}{\arabic{section}.\arabic{theorem}}

\newenvironment{remark}{\begin{protoremark}\rm}{\end{protoremark}}
\newenvironment{definition}{\begin{protodefinition}\rm}{\end{protodefinition}}
\renewenvironment{proof}{\par \trivlist
 \itemindent\parindent \item[\hskip\labelsep\sc Proof.]
 \ignorespaces}{\qed\endtrivlist}
\newenvironment{proofof}[1]{\par \trivlist
 \itemindent\parindent \item[\hskip\labelsep\sc Proof of #1.]
 \ignorespaces}{\qed\endtrivlist}

\newcommand\mnote[1]{}
\newcommand\bel[1]{\begin{equation}\label{#1}}
\newcommand\ee{\end{equation}}
\newcommand\lb[1]{\label{#1}\mnote{#1}}

\newcommand{\eps}{\varepsilon}
\newcommand{\GG}{\mathcal G}
\newcommand{\RR}{\mathcal R}
\newcommand{\QQ}{\mathcal Q}
\newcommand{\Z}{\mathbb Z}
\newcommand{\FF}{\mathbb F}
\newcommand{\HH}{{\mathcal H}}
\newcommand{\one}{\operatorname{\mathbf 1}}
\renewcommand{\mod}{\operatorname{mod}}
\newcommand{\Var}{\operatorname{Var}}

\newcommand{\GW}{\text{\sf GW}}
\newcommand{\ev}{\text{\bf E}}
\newcommand{\pr}{\text{\bf P}}
\newcommand{\dimum}{\dim_{\overline{M}} }
\newcommand{\dimh}{\operatorname{\dim_H}}
\newcommand{\Sym}{\operatorname{Sym}}
\newcommand{\Alt}{\operatorname{Alt}}

\newcommand{\reals}{{\mathbb R}}

\usepackage{natbib}
\begin{document}
\renewcommand{\thetheorem}{\arabic{theorem}}
\title{Dimension and randomness in groups  acting on \\ rooted trees}
\author{by Mikl\'os Ab\'ert\footnote {Department of Mathematics, University of Chicago, Chicago,
IL 60637, USA. {\sf abert@math.uchicago.edu}. Research
partially supported by OTKA grant T38059.}\;\ and B\'alint
Vir\'ag\footnote {Department of Mathematics, MIT, Cambridge, MA
02139, USA. {\sf balint@math.mit.edu}. Research partially
supported by NSF grant \#DMS-0206781.} }
\date{February 11, 2003} \maketitle
\begin{abstract}
We explore the structure of the $p$-adic automorphism group
$\Gamma$ of the infinite rooted regular tree. We determine the
asymptotic order of a typical element, answering an old
question of Tur\'an.

We initiate the study of a general dimension theory of groups
acting on rooted trees. We describe the relationship between
dimension and other properties of groups such as solvability,
existence of dense free subgroups and the normal subgroup
structure. We show that subgroups of $\Gamma$ generated by
three random elements are full-dimensional and that there exist
finitely generated subgroups of arbitrary dimension.
Specifically, our results solve an open problem of Shalev and
answer a question of Sidki.
\end{abstract}

\section{Introduction}

Let $T=T(d)$ denote the infinite rooted $d$-ary tree and let
$H\subseteq \mathrm{\ Sym}(d)$ be a permutation group. Let
$\Gamma (H)$ denote the infinite iterated wreath product of $H$
acting on $T$ with respect to $H$.
For example, $\Gamma(\mbox{Sym}(d))$ is the full automorphism group of $T(d)$%
. Let $\Gamma _{n}(H)$ denote the $n$-fold wreath product of $H$, acting on $%
T_{n}$, the $d $-ary tree of depth $n$.

The case $H=C_{p}$, the cyclic group of order $p$ is of
particular interest. The pro-$p$ group $\Gamma (p)=\Gamma
(C_{p})$ obtained this way is called the group of $p$-adic
automorphisms. The group $\Gamma _{n}(p)$ is called the
\emph{symmetric $p$-group} of depth $n$, as it can also be
obtained as the Sylow $p$-subgroup of the symmetric group
$\mathrm{Sym}(p^{n})$.

The first goal of this paper is to analyze the elements of
$\Gamma (H)$. We answer a question of Tur\'{a}n (see
\cite{psz83a}), which asks for an analogy of the famous theorem
of \cite{et65} about the distribution of orders of random
elements in $\mathrm{Sym}(n)$.

Let $p^{K_{n}}$ denote the order of a random element of the symmetric $p$%
-group $\Gamma _{n}(p)$ and let $\alpha _{p}$ be the solution
of the equation $\alpha (1-\alpha )^{1/\alpha-1 }=1-1/p$ in the
open unit interval.

\begin{theorem}
\lb{turanprob}We have $K_{n}/n\rightarrow \alpha _{p}$ in
probability.
\end{theorem}

A similar result holds for general $H$ (see Section
\ref{otre}). The proof of this theorem depends on a new
measure-preserving bijection between conjugacy classes of
random elements and Galton-Watson trees, bringing the theory of
stochastic processes to bear upon the subject. Theorem
\ref{turanprob} was conjectured by \cite{psz83a}; they proved
the upper bound and that the variance of $K_{n}$ remains
bounded as $n\rightarrow \infty$.
\cite{puchta2} showed
that the limit in Theorem \ref{turanprob} exists.
In a related recent paper,
\cite{evans} studies the random measure given by the
eigenvalues of the natural representation of $\Gamma _{n}(p)$.

The next goal is to understand the subgroup structure of
$\Gamma (p)$. Since every countably based pro-$p$ group can be
embedded into $\Gamma (p)$, we investigate subgroups according
to their Hausdorff dimension. We show that the Hausdorff
spectrum for the closure of subgroups generated by $3$
elements is the whole unit interval. More precisely,

\newcommand{\tetszdimthm}{
For each $\lambda \in \left[ 0,1\right] $ there exists a
subgroup $G\subseteq\Gamma (p)$ generated by $3$ elements such
that the closure of $G$ is $\lambda$-dimensional.}

\begin{theorem}
\lb{tetszdim} \tetszdimthm
\end{theorem}

This result leads to the solution of a problem of
\cite{MR2001d:20026}, who asked whether a finitely generated
pro-$p$ group can contain finitely generated subgroups of
irrational Hausdorff dimension.

In the context of pro-finite groups Hausdorff dimension was
introduced by \cite{abercrombie94}. It was deeply analyzed by
\cite{bs97} (see also \cite{MR2001d:20026} and
\cite{MR99d:17025}). They show that for closed subgroups,
Hausdorff dimension agrees with the lower Minkowski (box)
dimension. In the case of $ \Gamma(H)$, this can be computed
from the size of the congruence quotients $
G_{n}=G\Gamma^{(n)}/\Gamma ^{(n)}$, where $\Gamma ^{(n)}$
denotes the stabilizer of level $n$ in $\Gamma $. In fact,
$\dimh G$ is given by the $ \liminf $ of the \emph{density
sequence}
\bel{density}
 \gamma_{n}(G)={\log |G_{n}|/\log |\Gamma _{n}|}\text{.}
\end{equation}
For instance, the closure of Grigorchuk's group, a subgroup
of ${\Gamma (2)}$, has dimension $5/8$.

Theorem \ref{tetszdim} is obtained using probabilistic methods.
These methods give the following analogue of the theorem of
\cite{dixon69} saying that two random elements of $\Sym(n)$
generate $\Sym(n)$ or the alternating group $\Alt(n)$ with
probability tending to 1. Perhaps surprisingly, no
deterministic construction is known for large subgroups of
$\Gamma _{n}(p)$ with bounded number of generators as in the
following theorem.

\def\dixonthm{Let $\varepsilon >0$. Let $G_{n}$ be the subgroup
generated by three random elements of the symmetric $p$-group
$\Gamma _{n}(p) $. Then $\pr\left( \,|G_{n}|>|\Gamma
_{n}(p)|^{1-\varepsilon }\right) \rightarrow 1$ as
$n\rightarrow \infty $.}

\begin{theorem}\lb{dixon}
\dixonthm
\end{theorem}

For the infinite group $\Gamma (p)$, this takes the following
form: with probability 1, three random elements generate a
subgroup with 1-dimensional closure (see Theorem
\ref{3random}). We conjecture that the same holds for two
elements instead of three. Since three random elements generate
a spherically transitive subgroup with positive probability,
and $\Gamma (p)$ is not topologically finitely generated, the
above result shows the existence of a rich family of
topologically finitely generated transitive full-dimensional
subgroups of $\Gamma (p)$.

The analysis of randomly generated subgroups allows us to
answer a question of \cite{sidki01}. He asked whether the
binary adding machine, an element that acts on every level of
$T(2)$ as a full cycle, can be embedded into a free subgroup of
$\Gamma (2)$. Call a subgroup $G\subseteq \Gamma (H)$ strongly
free if it is free and every nontrivial element fixes only
finitely many vertices. For instance, both the trivial group
and the cyclic subgroup generated by the adding machine have
this property.

\def\strongfreethm{
Let $G\subseteq \Gamma (H)$ be a strongly free subgroup and let
$g\in \Gamma (H)$ be a random element. Then the group
$\left\langle G,g\right\rangle $ is strongly free with
probability 1. }
\begin{theorem}
\lb{strongfree} \strongfreethm
\end{theorem}

Next, we explore a general dimension theory of groups acting on
the rooted trees $T(p)$. The first result in this direction
shows that closed subgroups of $\Gamma (p)$ are `perfect' in
the sense of Hausdorff dimension. Recall that $G^{\prime }$
denotes the derived subgroup of the group $G$.

\def\perfektthm{Let $G\subseteq \Gamma (p)$ be a closed
subgroup. Then we have $\lim \left( \gamma _{n}(G)-\gamma
_{n}(G^{\prime })\right) =0$. In particular, $\dimh G=\dimh
G^{\prime }$.}\begin{theorem}\lb{perfekt} \perfektthm
\end{theorem}

As a corollary, solvable subgroups are zero-dimensional, and
positive-dimensional closed subgroups cannot be
\emph{abstractly} generated by countably many solvable
subgroups.

Unlike the infinite group $\Gamma (p)$, the finite symmetric $p$-groups $%
\Gamma _{n}(p)$ do have large Abelian subgroups of density
$1-1/p$. However, they cannot be glued to form a large subgroup
of $\Gamma (p)$. The following theorem explains the background
of this phenomenon.

\def\solsumthm{
Let $G\subseteq \Gamma (p)$ be a solvable subgroup with
solvable length $d$. Then we have $\sum_{n =0}^{\infty }\gamma
_{n }(G)\leq Cd$ where $C$ is a constant depending
 on $p$ only. }
\begin{theorem}
\lb{solsum}\solsumthm
\end{theorem}

One motivation to study dimension of subgroups in $\Gamma (p)$
comes from the theory of just infinite pro-$p$ groups, which
are regarded as the simple groups in the pro-$p$ category. A
group $G$ is just infinite if every proper quotient of $G$ is
finite. By the characterization of \cite{MR2002f:20044} (which
is based on \cite{MR43:338}), groups with this property fall
into two
classes: one is the so-called branch groups, which have a natural action on $%
T(p)$ (see Section \ref{large} for details). \cite{boston}
suggests a direct connection between Grigorchuk's classes and
Hausdorff dimension; he conjectures that a just infinite
pro-$p$ group is branch if and only if it has an embedding into
$\Gamma (p)$ with positive-dimensional image.

In the following we contrast known results about branch groups
with new results on 1-dimensional subgroups. We first show the
following.

\def\egyszeruthm{Let $G\subseteq \Gamma (p)$ be a spherically
transitive 1-dimensional closed subgroup. Then every nontrivial normal
subgroup of $G$ is 1-dimensional.}
\begin{theorem}
\lb{egyszeru}\egyszeruthm
\end{theorem}

It would be interesting to see whether a similar result holds
for spherically transitive positive-dimensional closed subgroups.

\cite{MR2001i:20060} conjectured that every just infinite
pro-$p$ branch group contains a nonabelian free pro-$p$
subgroup. An affirmative result is known if the group is not
virtually torsion-free (see \cite{MR2001i:20058}). Another
beautiful result in this vein is given by \cite{MR2002f:20037}
on groups satisfying the Golod-Shafarevich condition. In light
of Boston's conjecture it is natural to formulate the
following.

\def\propsubconj
{Let $G\subseteq \Gamma (p)$ be a positive-dimensional closed
subgroup. Then $G$ contains a nonabelian free pro-$p$ subgroup.
}
\begin{conjecture}
\lb{propsub} \propsubconj
\end{conjecture}

This might be more attackable than Wilson's original
conjecture. A step in this direction is the following.

\def\trapropthm{
Let $G\subseteq \Gamma (p)$ be a spherically transitive
1-dimensional closed subgroup. Then $G$ contains a nonabelian free
pro-$p$ subgroup.}
\begin{theorem}
\lb{traprop}\trapropthm
\end{theorem}

Another measure of largeness for profinite groups is whether
they contain dense free subgroups (see \cite{MR2002g:20114} and
\cite{MR2001h:22017}). \cite {MR2001i:20060} showed that just
infinite pro-$p$ branch groups contain dense free subgroups.
Let $d(G)$ denote the minimal number of topological generators
for $G$, which may be infinite.

\def\densefreethm{
Let $G\subseteq \Gamma (p)$ be a closed subgroup of dimension 1 and
let $k\geq d(G)$. Then $G$ contains a dense free subgroup of
rank $k$.}
\begin{theorem}
\lb{densefree}\densefreethm
\end{theorem}

In other words, the free group $F_{k}$ is residually $S$ where
$S$ denotes the set of congruence quotients of $G$. \bigskip

A common tool throughout the paper is the so-called \emph{orbit
tree} of a subgroup $G\subseteq \Gamma (H)$, the quotient graph
of $T$ modulo the orbits of $G$. Orbit trees reflect the
conjugacy relation (see \cite {MR2002k:20046}) and can be used
to describe the structure of Abelian subgroups in $\Gamma (p)$
(see Section \ref{small}). The key observation leading to
Theorem \ref{turanprob} is the fact that the orbit tree of a
random element is a Galton-Watson tree. A crucial step in
Theorems \ref{dixon} and \ref{strongfree} is that the orbit
tree of a randomly generated subgroup has finitely many rays
with probability 1, i.e., the closure has finitely many orbits
on the boundary of $T$ (Proposition \ref{kopaszfa}). Theorem
\ref{perfekt} also traces back to estimates on orbit trees.

Another tool used in the paper is word maps between the spaces
$\Gamma (p)^{k}$. A nontrivial word $w\in F_{k}$ can be thought
of as a map $\Gamma (H)^{k}\rightarrow \Gamma (H)$ via
evaluation. \cite{bhat} proved that random substitution into a
word map gives 1 with probability 0, or equivalently, the
kernel of the map has Haar measure zero. As a result, random
elements of $\Gamma (H)$ generate a free subgroup with
probability 1. The key result leading to Theorem
\ref{densefree} is the following generalization.

\def\kisdimkerthm{The kernel of a word map is not
full-dimensional in $\Gamma(H)^k$.}

\begin{theorem}
\lb{kisdimker}\kisdimkerthm
\end{theorem}

Another generalization leading to Theorem \ref{strongfree} is
that random evaluation of a word map yields an element of
$\Gamma (H)$ which fixes only finitely many vertices with
probability 1.

\bigskip

The structure of the paper is as follows.  Section \ref{otre}
contains the proof of Theorem \ref{turanprob}. It also
introduces notation and discusses conjugacy classes and random
elements in $\Gamma(H)$. Section \ref{otrs} discusses the orbit
structure of random subgroups. In Section \ref{wm} we explore
word maps and prove Theorems \ref{strongfree} and
\ref{kisdimker}. Sections \ref{sgh} and \ref{srs} introduce the
technical tools needed for proving Theorem \ref{dixon}, which
is done in Section \ref{drs}. Section \ref{small} discusses
small subgroups, and covers Theorems \ref{perfekt} and
\ref{solsum}. Section \ref{large} is about high-dimensional
subgroups, and contains the  proofs of Theorems \ref{egyszeru},
\ref{propsub}, \ref{traprop}, and \ref{densefree}.


\renewcommand{\thetheorem}{\arabic{section}.\arabic{theorem}}

\section{Elements}\lb{otre}

This section studies the statistical properties of elements of
$\Gamma(p)$. Random elements are described via the family tree
of a Galton-Watson branching process. This allows us to answer
Tur\'an's question (see Theorem \ref{turanprob}).

We first describe our notation for (iterated) {\bf wreath
products}. Let $(H_i,X_i)$, $1\le \ell \le n$ be a sequence of
permutation groups where $n$ might be infinite. Let
$d_i=|X_i|$. For $0\le \ell\le n$, define the {\bf level}
$\ell$ as
$$
\partial T_\ell:=X_1 \times \ldots \times
X_\ell,$$ and the tree $T_n$ which has vertex set $V(T_n)$
given by elements of the disjoint union of all levels. The
sequence $w\in V(T_n)$ is a child of $v\in V(T_n)$ (or $v$ is
the parent of $w$) if $w$ begins as $v$ and has one extra
element. The edge set $E$ is the collection of $\{v,w\}$ for
all such $v,w$.

Consider a map $g$ that assigns to each vertex $v\in
\partial T_\ell$ for $\ell<n$ an element of $H_{\ell-1}$. Let
$\Gamma$ be the set of all such maps. An element $g\in \Gamma$
acts on $T_n$ bijectively via the rule
 \bel{szorzet}
(x_1,\ldots, x_\ell)^g:=(x_1^{g(())},s_2^{g((x_1))}, \ldots,
x_\ell^{g((x_1,\ldots, x_{\ell-1}))}).
 \end{equation}
and $\Gamma=H_1\wr \ldots \wr H_n$ is easily checked to be a
subgroup of Aut$(T_n)$. This construction also works in the
case $n=\infty$.

If $n$ is finite, then the group structure of $\Gamma_n$ does
not depend on the permutation representation of the last group
$H_n$ on $\Sym(n)$, only on the abstract structure of $H_n$.
This defines $H \wr K$ where $K$ is an abstract group. When we
talk about $g\in H_1 \wr H_2$, the wreath product of two
groups, we will often  use the shorthand $g(x)$ for $g((x))$,
where $x\in X_1$. The rule (\ref{szorzet}) yields the
multiplication rule
$$
(gh)(x)=g(x)h(x^g).
$$
We will use the notation $\Gamma_n(H)$ for the wreath product
of $n$ copies of a permutation group $(H,X)$, and the shorthand
$\Gamma(H)=\Gamma_\infty(H)$. The most important cases are the
{\bf symmetric $p$-group} $\Gamma_n(C_p)$ and the full
automorphism group $\Gamma_n(\Sym(d))$ of $T_n$.

If $G$ is a group, and $S\subset G$, then $\langle S \rangle$
denotes the group {\bf abstractly generated} by $S$, that is
the set of elements that can be obtained from $S$ via repeated
group operations. If $G$ is a topological group, then we call
the closure $\overline{ \langle S \rangle}$ the group {\bf
topologically generated} by $S$.

Let $G$ be a group acting on a tree $T$.

\begin{definition}\label{otg} The {\bf orbit tree  of the group} $G$ is
the quotient graph $T_G$ of $T$ modulo the set of orbits of
$G$.
\end{definition}

It is easy to check that orbit trees are in fact trees.

For the purposes of Proposition \ref{conjugacy} the orbit trees
of {group elements} elements need to contain more information
about the local structure of automorphisms. Let $v$ be a vertex
of a tree $T=X^\infty$, and let $k$ denote the length of the
orbit of $v$ under an automorphism $g$. Then $g^{k}$ acts on
the child edges of $v$, which are labeled naturally by
elements of $X$. Let $\ell_g(e)\subseteq X$ denote the orbit of
the child edge $e$. For a set (e.g. orbit) $\mathcal E$ of
edges at the same level, let $\ell_g({\mathcal E})=\ell_g(e)$
for the lexicographically smallest $e\in {\mathcal E}$.

\begin{definition} The {\bf orbit tree of a tree automorphism $g$}
is the quotient graph  $T_g$  of $T$ modulo the set of orbits
of $g$, together with the labeling $\ell_g$ of the edges of
$T_g$.
\end{definition}

The orbit tree of $g\in \Gamma(H)$ is an $H$-labeled tree,
defined as follows.

An {\bf $H$-labeled tree} is a rooted tree together with an
edge-labeling $\ell$ so that for $v\in V$, $\ell$ is a
bijection between the child edges of $v$ and the cycles of some
$h_v\in H$.

Two $H$-labeled trees are called {\bf equivalent} if (1) there
exist a graph isomorphism $v \mapsto v',$ $e\mapsto e'$ between
them, and  (2)  for each $v$, there exists an $h_v\in H$ so
that $\ell(e')=\ell(e)^{h_v}$ for all child edges $e'$ of $v'$.

The following extension of a theorem of \cite{sushchansky} is
not difficult to check, and motivates the notion of orbit
trees.

\begin{proposition}\lb{conjugacy}
Two elements of $\Gamma(H)$ are conjugates if and only if their
orbit trees are equivalent.
\end{proposition}


In most important cases, equivalence is only a little more than
graph isomorphism. Say two labellings of a tree are equivalent
if they define equivalent $H$-labeled trees. In the case
$H=\Sym(d)$, two labellings are equivalent iff the length of
the edge label cycles agree. In this case, it suffices to label
edges by the cycle lengths. If $H=C_2$, then any two labellings
are equivalent. If $H=C_p$ for $p$ prime, then each vertex has
either $1$ or $p$ children, and labeling is the same as
cyclicly ordering the children whenever there are $p$ of them.

We now turn to the description of the orbit trees of random
elements. Our convention is that unless otherwise specified,
the term {\bf random elements} denotes independent random
elements chosen according to uniform (in the infinite case
Haar) measure.

\begin{definition} {\bf Labeled Galton-Watson (GW) trees.}
Let $\bar L$ be
the set of finite sequences with elements from a set of labels
$L$, and let $\nu$ be a probability distribution on $\bar L$.
We define the probability distribution $\GW(\nu)$ on infinite
trees with edges labeled by $L$ by the following inductive
construction. In the first generation, we have $1$ individual.
Given the individuals at level $n$, each of them has a sequence
of child edges picked from the distribution $\nu$
independently. The other endpoints of the child edges form
generation $n+1$.
\end{definition}

Let $(H,X)$ be a permutation group, and let $\nu_H$ be the
distribution of the sequence of orbits of a uniform random
element.

\begin{proposition}[Random elements and GW trees] \lb{haargw}
If $g\in \Gamma(H)$ is a Haar random element, then the
distribution of the orbit tree $T_g$ is $\GW(\nu_H)$.
\end{proposition}

\begin{proof}
We prove this by induction. The inductive hypothesis is that
the first $n$ levels of $T_g$ have the same distribution as the
first $n$ levels of a $\GW(\nu_H)$ tree.

Indeed, suppose this is true for $n$. Let $w=(v, \ldots
,v^{g^{k-1}})\in T_g$ be an orbit of an element on level $n$ of
$T$. Then the child edges of $w$ in $T_g$ correspond to, and
are labeled by orbits of $g^k(v)$. But
$g^k(v)=g(v)g(v^g)\ldots g(v^{g^{k-1}})$, and the factors on
the right hand side are chosen independently uniformly from
$H$, which implies that their product also has uniform
distribution.
\end{proof}

The simplest consequence of this is the following asymptotic
law. Let $r(H)$ denote the permutation rank of $H$, i.e. the
number of orbits of the action of $H$ on pairs of elements of
$X$ (e.g. $r(H)=2$ if $H$ is 2-transitive on $X$).

\begin{corollary}\lb{survival} Let the group $H$ be transitive on $X$, and let
$g_n\in \Gamma_n(H)$ be random.  Then as $n\rightarrow \infty$,
$$
n\pr(g_n\text{ fixes a vertex on level }n) \rightarrow
2/(r(H)-1).
$$
\end{corollary}

\begin{proof}
Let $N$ be the number of fixed points of a random element of
$H$ on $X$; then it is well known that $\ev N=1$, $\ev N^2=r$,
so $\Var\,N = r-1$.

A vertex at level $n$ of $T_{g_n}$ corresponds to a fixed point
if and only if all of its ancestor edges are labeled by a
fixed point. This happens if it lies in the subtree of
$T_{g_n}$ gained by removing all descendant subtrees which have
an ancestor edge not labeled by a fixed point. Such a subtree
has Galton-Watson distribution with branching structure given
by the number of fixed points $N$ of a random element in $H$.
Since $\ev N=1$, this is a critical tree, so by a theorem in
branching processes (see \cite{an}) the probability that it
survives until generation $n$ is asymptotic to $2/(n\Var N)$,
proving the corollary.
\end{proof}

The following corollary is immediate.

\begin{corollary}\lb{1nofix} Let $H$ be transitive on $X$, and
let ${g}\in \Gamma(H)$ be chosen according to Haar measure.
Then ${g}$ fixes only finitely many vertices of $T$ with
probability 1.
\end{corollary}

Let for an integer $n$, let $p(n)$ denote the highest power of
$p$ dividing $n$.
Let $h\in H$ be a uniform random element, and let $p$ be a
prime. Let $\mu_p(k)=\mu_{H,p}(k)$ denote the expected number
of orbits $v$ of $h$ with $p(|v|)=k$.
Let $\hat \mu_p$ be the generating function for
$\mu_p$:
$$
\hat \mu_p (z) = \sum_{k=0}^\infty \mu_p(k) z^k
$$
and note that the sum has only finitely many non-zero terms.

\begin{proposition}
Let ${g}_n\in \Gamma_n$ be a random element, let
 \bel{p-order} \alpha_p =
\min_{\lambda>0} \frac{\log \hat \mu_p(e^\lambda)}{\lambda}.
 \end{equation}
Then ${\ev\, p(|{g}_n|) / n} \rightarrow \alpha_p$ and there
exist $c_1,c_2>0$ so that for all $n,\,k$ we have
$$
\pr[\left| p(|{g}_n|) - \ev\, p(|{g}_n|)\right| \ge k ]\le c_1
e^{-c_2k},
$$
in particular, $\Var(p(|{g}_n|))$ remains bounded.
\end{proposition}

The bounded variance statement is somewhat surprising, as in
most limit theorems variance increases with $n$ (in most cases
it is on the order of $n$). The asymptotics of the order of a
typical element is immediate from this result.

\begin{corollary} As $n\rightarrow\infty$, we have
$$\log |{g}_n| /n \rightarrow \sum  \alpha_p \log p
$$
in probability, where the sum ranges over primes $p$ dividing
$|H|$.
\end{corollary}

As a special case, we get an answer to Tur\'an's question,
which we restate here. Let $p^{K_{n}}$ denote the order of a
random element of the symmetric $p$ -group $\Gamma _{n}(p)$ and
let $\alpha _{p}$ be the solution of the equation
 \bel{psz}
 \alpha (1-\alpha )^{1/\alpha-1 }=1-1/p
 \end{equation}  in $(0,1)$.

\newtheorem*{thm.tur}{Theorem \ref{turanprob}}
\begin{thm.tur}
We have $K_{n}/n\rightarrow \alpha _{p}$ in probability.
\end{thm.tur}

Tur\'an's goal was to find a natural analogue of the theorem of
the \cite{et65} about the asymptotics of the order of a random
element in $\Sym(n)$. Their paper started the area of
statistical group theory. Theorem \ref{turanprob} was
conjectured by \cite{psz83a}; they proved the upper bound and a
linear lower bound with a different constant. \cite{puchta2}
shows that the limit exists. He also studies uniformly
randomly chosen \emph{conjugacy classes} in
\cite{puchta}. In a related recent paper
motivated by random matrix theory, \cite{evans} studies the
random measure given by the eigenvalues of the natural
representation of $\Gamma_{n}(p)$.

\begin{proofof}{Theorem \ref{turanprob}}
 $H$ is the cyclic group of order $p$, so
 $
 \hat \mu_p(z) = 1 + z(p-1)/p.
 $
Let $f(\lambda)= \log \hat \mu_p (e^\lambda)$. Then
(\ref{p-order}) equals the minimum of $f(\lambda)/\lambda$, and
setting the derivative to $0$ we get
 \bel{eql}
 f(\lambda)/\lambda = f'(\lambda).
 \end{equation}
 At equality, this
expression gives $\alpha$. The right hand side is the
logarithmic derivative of $\hat \mu_p(e^\lambda)$, so it equals
$\alpha =(e^\lambda(p-1)/p) / (1+ e^\lambda(p-1)/p)$. Using this, we
easily express $\lambda$ and then $f(\lambda)$ from $\alpha$.
Substitution into (\ref{eql}) and algebraic manipulations give
$\ev K_n/n \rightarrow \alpha_p$, and convergence in
probability follows since the variance of $K_n$ is bounded.
\end{proofof}

For the proof of Proposition \ref{p-order}, we need to
understand the order of an element in terms of its orbit tree.
A simple inductive argument shows that if $v\in T_g$ is an
orbit, then its length is given by the product of the lengths
of labels on the simple path $\pi((),v)$ from the root $()$ to
$v$ in $T_g$:
 \bel{order}
|v| = \prod_{e\in \pi((),v)} |\ell(e)|,
 \end{equation}
 and therefore
\bel{primediv} p(|v|)= \sum_{e\in \pi((),v)} p(|\ell(e)|).
 \end{equation}

Proposition \ref{ol} below shows that all properties about the
length of orbits of $\Gamma_n$ can be understood in terms of a
branching random walk.

Let $\{p_i\}_{1\le i \le d}$ be a sequence of primes. Consider
an orbit $o$ of an action group $H$ on a set $X$, and let
ex$(o,X)$ denote the point in $\Z^d$ whose $i$th coordinate is
the exponent of the highest power of $p_i$ dividing $|o|$. Let
ex$(H,X)$ denote the multiset $
 \{\text{ex}(o)\ | \ o\text{\ orbit of }H\}.$

\begin{definition} Let $\nu$ be a probability distribution
on the space $S$ of finite multi-sets of elements of $\mathbb
Z^d$. A {\bf branching random walk} is a probability measure on
infinite sequences $\{\Xi_n;\, n\ge 0\}$ with elements from
$S$. It is constructed recursively as follows. At time $n=0$ we
have $\Xi_0=\{0\}$, the set containing the origin.

If $\Xi_n$ is already constructed, then each $x\in \Xi_n$ has
``offspring'' $x+Y_1, \ldots, x+Y_N$. Here $(Y_1, \ldots, Y_N)$
is picked from $\nu$ independently from the past and from the
offspring of other individuals at time $n$. $\Xi_{n+1}$ is the
multi-set union of all the offspring of $x\in \Xi_n$.
\end{definition}

We call the measure
$$
\mu(x)= \ev\left( \sum_{i=1}^N \one(Y_i=x) \right)
$$
the {\bf occupation measure} of the BRW. Let $\hat \mu$
denote the generating function for the measure $\mu$.

In our case, the relative offspring positions are given by the
sequence
$$\{\,(p_1(|o|),\ldots, p_d(|o|))\, :\, o\text{ orbit of }h\,\},$$
where $h$ is a uniform random
element of $H$.

\begin{proposition}[Orbit lengths and BRW]\lb{ol}
Let ${g} \in \Gamma(H)$ be a Haar random element, and let
$\Xi_n=\text{ex}(\langle{g}\rangle,\partial T_n)$. Then
$\{\Xi_n;\, n \ge 0\}$ is a branching random walk on $\Z^d$.
\end{proposition}

\begin{proof} This is immediate from the formula
(\ref{primediv}) and the proof of Proposition \ref{haargw}.
\end{proof}

Since $p(|g|)$ equals the maximum of $p(|o|)$ for the orbits of
$g$, we need to understand the position of the highest element of
$\Xi_n$. The corresponding questions for the branching random
walk have been answered by \cite{hammersley}, \cite{biggins},
and \cite{dh91}. The part of their results that is relevant to
our setting is stated in the following theorem.

\begin{theorem}\lb{brw}
Let $X_n$ denote the position of the greatest individual of a
1-dimensional a BRW whose occupation measure $\mu$ has finite
support. Let
$$
\alpha= \inf_{\lambda>0}\,\frac{\log
\hat\mu(e^{\lambda})}{\lambda}
$$
Then
$$
\ev X_n/n \rightarrow \alpha,
$$
and there exist $c_1,c_2>0$ so that for all $n$,
$$
\pr[|X_n-\ev X_n| \ge k] \le c_1 e^{-c_2n}.
$$
\end{theorem}

Here we present a sketch of the proof of the first statement in
a simple case.

\begin{proof}
Let $\mu^{*n}(k)$ denote the expected number of individuals at
position $k$ at time $n$. If we find the highest $k$ for which
this quantity is about $1$, that means that on average there
will be $1$ individual at position $k$, and less than one
individuals at positions higher than $k$.  Intuitively, this
$k$ seems to be a good guess for the position of the highest
individual. There are technical arguments to make this
intuition rigorous.

The next step is to solve $\mu^{*n}(k)=k$. By considering one
step in the branching random walk one discovers the convolution
formula
$$
\mu^{*n}(k)= \sum_{j=0}^k \mu^{*n-1}(j)\mu(k-j)
$$
but convolution turns into product for generating functions:
$\widehat {\mu^{*n}} =  \hat \mu^n$. For example, for the
symmetric $p$-group
$$ \hat \mu(z)^n = (1+z(p-1)/p)^n,
$$
and we are left with finding the power $k$ of $z$ for which the
coefficient in the above expression is about $1$. This gives
the approximate equation
 \bel{psz2}
  \binom{n}{k} = \left(\frac{p}{p-1}\right)^k, \ \ \ \ \alpha=k/n
 \end{equation}
and by the Stirling formula the asymptotic version of this
equation is (\ref{psz}). For the general case, the theory of
large deviations is a standard tool for approximating
$\mu^{*n}(\alpha n)$ for large $n$.
\end{proof}

\section{Orbit trees of random subgroups}
\lb{otrs}

Just as for single elements, there is a natural graph structure
on the orbits of a subgroup $G\subseteq \Gamma(H)$, see
Definition \ref{otg}. In this section, we prove that the orbit
tree of random subgroup (i.e. subgroups generated by random
elements) is a multi-type Galton-Watson tree. As a result, we
show that the closure of a random subgroup has finitely many
orbits on $\partial T$ with probability 1, and compute the
probability that it is transitive.

\begin{definition} {\bf Multi-type GW trees}. Let $\Upsilon$ be
a finite or countable set called {\bf types}. Let $\bar
\Upsilon$ be a set of finite sequences of $\Upsilon$. For each
$y \in \Upsilon$ let $\mu_y$ be a probability measure on $\bar
\Upsilon$. A multi-type Galton-Watson tree $\GW(\mu_\cdot)$ is
a random tree constructed as follows. The first generation
consists of an individual is of a given type $y\in \Upsilon$.
Individuals in generation $n$ have children according to the
measure $\mu_y$, where $y$ is the type of the individual. The
children of a generation are picked independently and they form
the next generation.
\end{definition}

Let $j\ge 1$, $k\ge 0$, and $(H,X)$ a permutation group.
Consider the action of the subgroup generated by $(j-1)k+1$
random elements of $H$ on $X$. Consider the list of the length
of the orbits, and replace each element $\ell$ of this list by
$\ell k$. The distribution of the new list yields a probability
measure $\mu_{H,j,k}$ on finite sequences of integers.

\begin{proposition}\lb{multi} Let $G$ be the subgroup generated by $j$
random elements of $\Gamma(H)$. Then $T_G$ has multi-type
Galton-Watson distribution $\GW(\mu_{H,j,\cdot})$, where the
types correspond to the length of the orbits.
\end{proposition}

First we introduce some notation. Let $G=H \wr K$ a permutation
group, let  $G_x$ denote the stabilizer of $x$, $G_{x+}$ the
subgroup of elements $g\in G_x$ for which $g(x)=id$. The group
$G_x/G_{x+}$ ``ignores'' the action everywhere except at $x$;
it is naturally isomorphic to $K$.

The following observation is immediate.

\begin{fact}\lb{stab}
Let $G\subseteq \Gamma(H)$,  and let $T_{G,v}$ be the
descendant subtree of the $T_G$ at the orbit $w$ of $v\in T$.
Then $T_{G,v}$ is isomorphic to $T_{G_v/G_{v+}}$ with the
labels multiplied by the length of $w$.
\end{fact}

\begin{lemma}\lb{pi1}
Let $(H,X)$ be a permutation group and let $K$ be group.
Consider the subgroup $G$ generated by $j$ Haar random elements
of $H \wr K$. Conditioned on the orbit $V$ of $x$, $G_x/G_{x+}$
has the same distribution as a subgroup generated by
$|V|(j-1)+1$ random elements of $G$.
\end{lemma}

\begin{proof}

Let $g_i$ denote the randomly chosen generators of $H$.
Consider the directed Schreier multi-graph $(V,E)$ of the
action of $G$ on $V$. For each element $v\in V$ this graph has
out-edges $(v,i)$ from $v$ to $v^{g_i}$. Let $\pi_1$ denote its
fundamental group of the graph. It consists of equivalence
classes of cycles starting from $x$. The cycles may contain an
edge in either direction. The equivalence relation is generated
by adding or removing immediately retraced steps.

It is well-known that for a connected graph, $\pi_1$ is a free
group generated by $|E|-|V|+1$ elements. Here is one way to
specify a set of generators. Take a spanning tree of $(V,E)$.
For each edge $e$ in the set $E_0$ of edges  outside the
spanning tree, the path $w_e$ is follows the path in the tree
from $x$ to $e_-$, moves along $e$ and then back in the tree to
$x$. Clearly, $|E_0|=|E|-|V|+1$.

Each path from $x$ is described by a word in the generators
$g_i$ and their inverses. If two paths are equivalent, then the
corresponding words $w$ have evaluate at $x$ to the same
element of $K$. Moreover, a word is in $G_x$ if and only if it
corresponds to a cycle.

It follows that $G_x$ is generated by the words corresponding
to generators of $\pi_1$. It suffices to show that for the
generators discussed above, $\{w_e/G_{x+}\}_ { e\in E_0}$ are
independent random elements of $G_x/G_{x+}$.

For an edge $e=(v,i)\in E$, let $g(e)=g_i(v)$. Respectively,
let $g(e)=g_i(v)^{-1}$ for the reversed edge. If $w=e_1 \ldots
e_k$ is a cycle starting from $x$, then
$$
w/G_{x+} = g(e_1) \ldots g(e_k).
$$
Now fix the values of $g(e)$ corresponding to $e\in E \setminus
E_0$. Then for $e\in E_0$, $w_e/G_{x+}$ equals $g(e)$
multiplied on the left and on the right by fixed elements.
Since the $\{g(e)\}_{ e\in E_0}$ are uniform random and
independent, so are the $\{w_e/G_{x+}\}_{ e\in E_0}$, as
required.
\end{proof}

\begin{corollary}\lb{pi1c}
For $v\in T$ of orbit length $k$ the subgroup $G_v/G_{v+}$ has
the same distribution as $G_{(j-1)k+1}$.
\end{corollary}

The following fact is also immediate.

\begin{fact}\lb{ind}
Let $v_1, \ldots v_k$ be vertices of $T$. Let $A$ denote the event
that none of the vertices $w(v_i)$ of $T_{H}$
equals another or
an ancestor of another. Conditioned on $A$, the subgroups
$G_{v_i}/G_{v_i+}$ are independent.
\end{fact}

\begin{proofof}{Proposition \ref{multi}}
Proposition \ref{multi} easily follows from Corollary
\ref{pi1c}, Fact \ref{stab}, and Fact \ref{ind}.
\end{proofof}

Let $q(j)$ denote the probability that $j$ random elements
generate a transitive subgroup $G$ of $H$. Looking at the
probability that $T_G$ does not branch we get the following

\begin{corollary} \lb{ptransc}
 \bel{ptrans}
 \pr(G\mbox{ is transitive on level }n)=
 \prod_{\ell=0}^{n-1} q\left(1+(j-1)|X|^\ell\right).
 \end{equation}
\end{corollary}
A subgroup of $\Gamma(H)$ is called {\bf spherically
transitive} if it acts transitively on every level of the tree.
As an example of Corollary \ref{ptransc}, consider two random
elements of the automorphism group of the binary tree. In this
case $j=2$, $H=C_2$, we have $q(n)=1-2^{-n}$, and
$$
\pr(G \mbox{ is spherically transitive}) =
\prod_{\ell=0}^\infty (1-2^{-2^\ell-1})\sim 0.63.
$$
In the general case, we have
\begin{corollary}\lb{trans}  Let $j>1$, and assume that $H$ has a transitive subgroup
generated by $j$  elements. Then with positive probability, $j$
independent Haar random elements of $\Gamma(H)$ generate a
spherically transitive subgroup. Moreover, this probability
tends to $1$ as $j \rightarrow \infty$.
\end{corollary}

\begin{proof}
Clearly, $1-q(nj) \le (1-q_j)^n$. This means that
$1-q(1+(j-1)|X|^\ell)$ is summable, so for $n=\infty$ the
product (\ref{ptrans}) is positive.
\end{proof}

\begin{proposition}\lb{kopaszfa}
Let $G$ be the subgroup generated by two independent Haar
random elements of $\Gamma(H)$. Then with probability 1, $T_G$
has only finitely many rays. Equivalently, the number of orbits
at level $\ell$ remains bounded as $\ell \rightarrow \infty$.
\end{proposition}

\begin{proof}
Consider the vertices of $T_G$ that are orbits of length 1 in
$T$. These induce a subtree with a subcritical GW distribution
(the number of fixed points of just one random element has
expectation 1 by P\'olya theory).

This implies that with probability 1  (1) there exists finitely
many fixed vertices of $T$, (2) there exists finitely many
vertices of $T$ with orbits of length at most $k$ for each $k$.

On the basis of Corollary \ref{trans} let $k$ be so large that
the chance that $(|X|-1)k+1$ random elements generate a
transitive subgroup of $\Gamma(H)$ is at least $1-1/|X|$.

Now consider the tree built to a level where each orbit has at
least $k$ elements. Then each such orbit will stop splitting
forever with probability at least $1-1/|X|$, and if it does
not, then it splits into at most $|X|$ offspring with
probability at most $1/|X|$. Thus the orbit process is
dominated by the subcritical GW tree with $0$ or $|X|$
offspring with probability at least $1-1/|X|$ and at most
$1/|X|$, respectively. Hence all orbits stop splitting with
probability 1.
\end{proof}

\section{Word maps}\lb{wm}

Let $G$ be a group with finite Haar measure. A {\bf randomly
evaluated word} in a group $G$ is a $G$-valued random variable
defined by a word $w$ in the free group $F_k$ by substituting
$k$ independent Haar random elements of $G$ into the generators
of $F_k$.

Randomly evaluated words are perhaps the simplest possible maps
$G^k \rightarrow G$. \cite{bhat} showed that with probability
1, a nontrivial word randomly evaluated in $\Gamma(H)$ does not
give the identity; moreover, it does not stabilize a fixed ray
in the tree. In this section we strengthen this result in two
directions, and give a sufficient condition for a list of
randomly evaluated words to have independent uniform
distribution.

Our first result is key in answering a question of
\cite{sidki01}.

\begin{proposition}\lb{wordfix}
A nontrivial randomly evaluated word in $\Gamma(H)$ fixes only
finitely many vertices of $T$ with probability 1.
\end{proposition}

\begin{proof}
Let $w= w_1 \ldots w_\ell$, where each $w_i$ is one of the
random generators or its inverse, and $w$ is in reduced form.
We use induction on the length $\ell$ of $w$. The $\ell=1$ case
is exactly Corollary \ref{1nofix}.

Let $\overline w$ denote the random evaluation of $w$. Let
$v\in T$ be a vertex. Let $w'_i=w_1\ldots w_i$,
$v_i=v^{\overline w'_i}$. Make the assumption $A$ that
$v^{\overline w}=v_\ell=v$, but $v_i \not=v_j$ for $0\le
i<j<\ell$. We use the notation $g(v)$ for the action of $g\in
\Gamma(H)$ on the descendant subtree of $v$. Then
$$
{\overline w}(v)= \prod_{i=1}^{\ell} {\overline
w}'_i(v^{{\overline w}'_{i-1}}).
$$
Given the event $A$, the factors in the product are independent
Haar random elements of $\Gamma(H)$,  and therefore ${\overline
w}(v)$ is also a Haar random element of $\Gamma(H)$. This means
that with probability one, by Corollary \ref{1nofix}, only
finitely many of the descendants of $v$ are fixed by
${\overline w}$.

Now we are ready to prove the proposition. By the inductive
hypothesis, there are only finitely many vertices $v$ which are
fixed by the evaluation of any of the proper sub-words
$w_i\ldots w_j$, $i<j$ of $w$. Let $L$ be the greatest level at
which there is such a vertex. Then at level $L+1$ event $A$
holds for all vertices $v$ that are fixed by ${\overline w}$.
In particular, all vertices at this level have finitely many
offspring that are fixed by ${\overline w}$, as required.
\end{proof}

This proof generalizes to the infinite wreath product of
arbitrary (not necessarily identical) finite transitive
permutation groups. From probabilistic point of view, this
property is interesting because it describes a critical
phenomenon in the statistical physics sense. One manifestation
of this is that the probability that a word is satisfied on
levels $1$ to  $n$  decays polynomially in $n$, as opposed to
exponentially (see Corollary \ref{survival}).

The following corollary is immediate.

\begin{corollary}\lb{freefree} %
Let $G$ be a subgroup abstractly generated by $k$
random elements of $\Gamma(H)$. Then $G$ is a free group acting
freely on $\partial T$.
\end{corollary}

In fact, we have that all elements of the countable group $G$
have only finitely many fixed vertices in $T$. We call a
countable free subgroup $G$ of $\Gamma(H)$ with this property
{\bf strongly free}. We are ready to prove Theorem
\ref{strongfree}, which we restate here.

\newtheorem*{strongfree.thm}{Theorem \ref{strongfree}}
\begin{strongfree.thm}
\strongfreethm
\end{strongfree.thm}

\begin{proof}
The proof of Proposition \ref{wordfix} did not use the full
power of randomness of all generators; it was sufficient to
have one random generator.

In fact, we have the following corollary to the proof of
Proposition \ref{wordfix}. Let $g_1$ be random and $g_{2},
\ldots g_k$ fixed elements of $\Gamma(H)$. Let $w$ be a reduced
word in the $g_i$, containing $g_1$ or its inverse. If the
evaluation all proper sub-words of $w$ fixes finitely many
vertices of $T$ with probability 1, then the same is true for
$\overline w$.

Theorem \ref{strongfree} follows.
\end{proof}

The {\bf adding machine} $s$ is an element of $\Gamma(p)$ that
acts as a full cycle on $\Gamma(p)$, or equivalently, one that
has orbit tree given by a single ray. More precisely, it is a
specific such element, defined by the following identity. For a
vertex $v=(v_1, \ldots, v_n)$ at level $n$ let $N(v) = v_1 +
v_2 p+ \ldots +v_n p^{n-1}$. Then for every $n$ and vertex $v$
at level $n$ we have
 $$
 \qquad N(v^s)=N(v)+1 \qquad (\mod p^{n}).
 $$
The adding machine $s$ is one of the simplest examples of
elements of $\Gamma(p)$ that are given by finite automata.
Groups generated by such elements include Grigorchuk's group,
and have a rich theory, see \cite{gns00}. \cite{sidki00} showed
that $s$ and another element given by a finite automaton cannot
generate a free subgroup of $\Gamma(p)$.  He asked
(\cite{sidki01}) whether the adding machine and {\it any} other
element can generate a free group. Since the adding machine
generates a strongly free cyclic subgroup, we have

\begin{corollary}[Answer to a question of Sidki] The adding
machine and a random element abstractly generate a free
subgroup of $\Gamma(p)$ with probability 1.
\end{corollary}

Simple maps $\reals^k \rightarrow \reals$ tend to have the
property that the pre-images of points are at most $k-1$
dimensional.  A word in $k$ letters evaluated in $\Gamma(H)$ is
a simple map $w:\Gamma(H)^k \rightarrow \Gamma(H)$. It is
natural to ask whether the analogous statement is true here.
The answer is no, for $k=1$ a counterexample is the map $x
\mapsto x^{p^k}$ in $\Gamma(H)$.  If $g\in \Gamma(H)$ acts as a
full cycle $(v_1, \ldots, v_{p^\ell})$ on level $\ell$ of the
tree, and the action $g(v_i)$ on the descendant subtrees
satisfy $g(v_1)\ldots g(v_{p^\ell})=1$, then $g$ is in the
kernel of this map. The set of such elements has dimension
$1-p^{-\ell}$.

For general $k$, we first need to define dimension. Consider
the usual $L_1$-metric on the product space $\Gamma(H)^k$,
which is built from the metric \eqref{metric} on $\Gamma(H)$.
This automatically gives rise to a notion of dimension in the
usual way. If $S\subseteq \Gamma(H)^k$, then it is easy to
check that upper Minkowski dimension, which dominates Hausdorff
dimension, satisfies \bel{dimum}
 \dimum S = \limsup_{n\rightarrow \infty} \left( \log
| S/(\Gamma^{(n)})^k|\,  /\log |\Gamma_n|\right).
\end{equation}

The 1-dimensional example above easily generalizes to show that
the  kernel of a $k$-letter word map can have Hausdorff
dimension arbitrarily close to $k$. However, the kernel of a
map cannot be full-dimensional, as the following restatement of
Theorem \ref{kisdimker} shows.

\begin{theorem} Let $w$ be a nontrivial word in the
letters $g_1,\ldots, g_k$, and let
$$\mathcal K=\{(g_1,\ldots, g_k)\in \Gamma(H)^k\,:\, w=1\}.$$
Then $\dimh(\mathcal K)\le \dimum(\mathcal K)<k$.
\end{theorem}

\begin{proof} The first inequality holds in every metric space.
We first make the assumption (i) that there exists elements in
$H$ for which $w\not=1$. Let $\ell$ denote the length of the
word (i.e., the sum of the absolute values of the exponents),
and let $b=|X|$. Fix the action $A_n=(g_i/\Gamma^{(n-1)};\;1\le
i \le k)$ of the variables on $T_{n-1}$,  let $s(A_n)$ denote
the number of liftings of these elements to $\Gamma_n$ for
which $w=1$. Let $h=|H|$, and let $t_n = h^{kb^n}$ denote the
total number of liftings. Let $s_n=\max_{A_n} s(A_n)$. Then we
have
$$
|\mathcal K/(\Gamma^{(n)})^k| \le s_1\ldots s_n,
$$
and therefore by \eqref{dimum}
$$ \dimum \mathcal K \le k \limsup \frac{\log
(s_1\ldots s_n)}{\log (t_1 \ldots t_n)} \le k \limsup
\frac{\log s_n}{\log t_n}.
$$
The last inequality is simple arithmetic.

Since the action on $T_{n-1}$ is fixed, for $v\in
\partial T_{n-1}$ we have that $w(v)$ is a word in the letters
$g_i(u),\, u\in \partial T_{n-1}$. More precisely, we get
$w(v)$ by replacing each occurrence of $g_i$ in $w$ by some
$g_i(u)$, where the $u$ depend on the action $A_n$ and need not
be the same for different occurrences of $g_i$. Together with
our assumption (i), this implies that there is a counterexample
to $w(v)$ in $H$.

We claim that there exists a large subset $S\subseteq
\partial T_{n-1}$ so that if $u,v\in S$, then the words $w(u),\,w(v)$
have no letter in common. Indeed, each word uses at most $\ell$
letters, and conversely, each letter appears in at most $\ell$
words (more precisely, the number of words using $g_i(u)$ is at
most the total number of times $g_i$ is used in $w$). This
means that the graph in which $u,v$ are neighbors iff $w(u)$,
$w(v)$ share a letter  has maximal degree at most
$d=\ell(\ell-1)$ (ignoring loops). Such a graph has an
independent set $S$ of size at least $\lceil |\partial
T_{n-1}|/(d+1)\rceil$ (using the greedy algorithm).

Let $v\in S$ and consider the set of letters $L$ that
contribute to $w(v)$.  The number of assignments of values in
$H$ to the letters $L$ is $h^{|L|}$; the number of assignments
for which $w(v)=1$ is at most $h^{|L|}-1$. Since these set of
letters are disjoint for different $v$, we have that
$$
s(A_n) \le t_n ((h^\ell-1)/h^\ell)^{|S|}
$$
and therefore with $\eps= \log (h^\ell/(h^\ell-1))$, we get
$$
\log s_n \le \log t_n - \eps b^n / (d+1) \le (1-\eps')\log t_n,
$$
where $\eps'>0$ does not depend on $n$.  This implies the claim
of the proposition when assumption (i) holds. Otherwise, for
some finite $j$ the group $H'=\Gamma_j(H)$ has a
counterexample, since $\Gamma(H)$ contains a free subgroup of
rank $k$ (see \cite{bhat} and Corollary \ref{freefree} here).
Since the natural isomorphism between $\Gamma(H)^k$ and
$\Gamma(H')^k$ preserves full-dimensionality, the proposition
applied to $H'$ concludes the proof.
\end{proof}
\begin{remark}
We have shown that the pre-image of $1\in \Gamma(H)$ is not
full-dimensional. It is easy to modify this proof to get that
the pre-image of any point $g\in \Gamma(H)$ has dimension at
most $1-\eps'(w)$, where $\eps'(w)$ does not depend on $g$.
\end{remark}

\begin{corollary}\lb{1dfree} Let $G\subseteq \Gamma(H)$ be a subgroup of
Hausdorff dimension $1$. Then $k$ random elements generate a
free subgroup with probability 1.
\end{corollary}

\begin{proof}
It suffices to show that for each word $w$ in $k$ letters, the
probability that $\overline w=1$ in $G$ is zero. This
probability equals the measure of $\mathcal K_w$ in $G^k$. But
$G^k$ is a full-dimensional subgroup of $\Gamma(H)^k$, so every
subset of $G^k$ with positive measure has full dimension in
$\Gamma(H)^k$.
\end{proof}

Let $G$ be a finite $p$-group. Our study of random generation
requires a sufficient condition for a set of randomly evaluated
words to have independent uniform distribution on $G$. For a
word $w$ in the variables $g_i$ define the {\bf exponent sum
vector} of $w$ as the vector whose coordinate $i$ is the sum of
the exponents of $g_i$ in $w$. Let $\overline w$ denote the
evaluation of the word $w$. We have not been able to locate the
proof of the following simple lemma in the literature.

\begin{lemma}[Linear and probabilistic independence]\lb{lip}
Let $w_1, \ldots, w_k$ be words
in the letters $g_1, \ldots, g_{n+m}$, so that the exponent sum
vectors restricted to the first $n$ coordinates are linearly
independent mod $p$. Fix $g_{n+1}, \ldots g_{n+m} \in G$. Then
the substitution map $G^n \rightarrow G^k$,
$$
(g_1, \ldots, g_n) \mapsto (\overline w_1, \ldots, \overline
w_k)
$$
covers $G^k$ evenly, i.e. it is $|G|^{n-k}$-to-one and onto.
\end{lemma}

A special case of Lemma \ref{lip} says that if $w$ is a word in
which the exponent sum of some letter is relatively prime to
the order of a finite $p$-group $G$, then the evaluation map
covers $G$ evenly. We state without proof that this also holds
for nilpotent groups but not for arbitrary finite groups.

\begin{proof} Assume that the $g_1, \ldots, g_n$ are chosen
independent uniformly and at random (i.u.).  It suffices to
show that the evaluations $\overline w_i(g)$ are i.u. We prove
this by induction; it is clearly true for $G=1$. Let
$Z\vartriangleleft G$ be a subgroup of order $p$ of the center
of $G$. Let $R$ be a set of coset representatives of $Z$, and
let $r: G \rightarrow R$ be a coset representative map.

Write $g_i=z_ir(g_i)$. Then $z_i\in Z$, $r_i\in R$ are i.u. for
$i\le n$. Because the $z_i$ are in the center, we have
$$
{\overline w}_i(g)={\overline w}_i(z){\overline
w}_i(r(g))={\overline w}_i(z) b_i r({\overline w}_i(g))
$$
Where $b_i\in Z$ and the random vector $b$ depends on $r(g)$,
but not on $z$. The assumption implies that ${\overline
w}_i(z)$ are i.u., and therefore given $b$, ${\overline w}_i(z)
b_i$ are i.u. Thus ${\overline w}_i$ are products of i.u.
elements of $Z$ and i.u. coset representatives, and these $2n$
random variables are also jointly independent. The claim
follows.
\end{proof}

By considering finite quotients, we see that the assumption of
Lemma \ref{lip} implies that for any pro-$p$ group $G$ the
substitution map is a measure-preserving surjection.

\section{Schreier graphs and homology}\lb{sgh}

The goal of this section is to establish some tools for proving
Theorem \ref{dixon} about random generation.

Let $K$ be a $p$-group acting on a set $X$, and let
$S=\{g_1,\ldots ,g_\ell\}$ be a multi-set of elements of $K$.
The {\bf Schreier graph} $\GG=\GG(X,S)$ is a directed
multi-graph whose vertex set is $X$, and the edge set is
$X\times S$, where the edge $(v,g)$ connects $v$ to $v^g$.

Let  $\tau \GG$ denote the vector space of {\bf 1-chains}, that
is the set of abstract linear combinations of the edges of
$\GG$ with coefficients in $\FF_p$.

If $w$ is a word in $S$, i.e. an element of the free group
indexed by $S$, then for each $v\in X$, $w(v)$ can be thought
of as a path $(v,v^{\overline w_1},v^{ \overline w_n})$, where
the $w_1, \ldots, w_n$ are the initial sub-words of $w$ and
$\overline w$ denotes evaluation in $K$.

For a path $\pi$, let $\tau \pi\in \tau \GG$ denote the {\bf
1-chain of $\pi$} where each edge appears with a coefficient
given by the number of times it is used in $\pi$.

Recall that the {\bf exponent sum vector} of $w$ is the vector
whose coordinate $i$ is the sum of the exponents of $g_i$ in
$w$.

\begin{lemma}[Independent 1-chains]\lb{homo} Let $W$ be a
set of words in whose exponent sum vectors
are linearly independent. Then $\{\tau w(v)\ :\ w\in W, v\in
X\}$ are linearly independent.

 More generally, let $W$ be a set of
words in $S$ whose exponent sum vectors restricted to
$S'\subseteq S$ are linearly independent. Then $\{\tau w(v)\ :\
w\in W, v\in X\}$ restricted to $S'$-edges are linearly
independent.

\end{lemma}

\begin{proof} We prove the more general version.  Let $S=\{g_1,\ldots,
g_\ell\}$, let  $n\le \ell$,  and let $S'=\{g_1,\ldots, g_n\}$.
Without loss of generality we may extend the set $W$ so that
$|W|=n$. Consider the wreath product $D=C_p \wr K$, and
consider the evaluation map
\begin{eqnarray*}
D^n&\rightarrow& D^n \\
(d_1,\ldots, d_n)&\mapsto &(\overline w_1,\ldots, \overline
w_n)
\end{eqnarray*}
where for $i>n$ we take the $i$th letter to be constant
$d_i:=(0,g_i)$. This map is a bijection because of the linear
independence assumption and Lemma \ref{lip}. Now restrict the
mapping to all possible extensions of
 $(g_1,\ldots, g_n)\in K^n$ given by the fixed generators so that
$(h_1,\ldots, h_\ell)\mapsto (h'_1,\ldots,h'_n)$ where
$h_i=(\nu_i,g_i),\, h'_i=(\nu_i',g'_i)$. Then the $g_i'$ are
determined by the $g_i$, and the map $f:\nu \mapsto \nu'$ is an
isomorphism of vector spaces $\FF_p^{|X|n} \rightarrow
\FF_p^{|X|n}$. If the reduced form $w\in W$ is $a_1\cdots a_m$,
then for $v\in X$ the evaluation in $D$ satisfies
 \bel{szorzet2} {\overline w}(v) =
{\overline a}_1(v){\overline a}_2(v^{\overline a_1}) \cdots
{\overline a}_m(v^{{\overline a}_1\cdots {\overline a}_{m-1}}).
 \end{equation}
Here $(v)$ means coordinate in the wreath product. Since the
edges of the subgraph $\GG(X,S')$ are of the form
$\{(v,g_i):v\in X, 1\le i \le n\} $, which agree with the
coordinates of $\nu$, we may think of $\nu$ as a 1-chain for
$\GG(X,S')$. Then (\ref{szorzet2}) implies that for $w\in W$,
$v\in X$ we have
$$ f(\nu)(w,v)= \overline w(v) = \tau' w(v) \cdot \nu,
$$
where $\tau' w(v)$ is the restriction of $\tau w(v)$ to
$S'$-edges, and $\cdot$ is the natural scalar product of
1-chains in $\GG(X,S')$. Since $f$ is bijective, the statement
of the lemma follows.
\end{proof}

The subspace of $\tau \GG$  generated by $\tau \pi$ for cycles
$\pi$ is called the (first) {\bf homology group} $\HH\GG$ of
$\GG$. For a word $w$ and $v\in X$ let $w(\circ v)$ denote the
path $(w^\ell)(v)$, where $\ell$ is the smallest positive power
so that $\overline w^\ell$ fixes $v$. Note that the path
$w(\circ v)$ is a cycle, and so $\tau w(\circ v)\in \HH\GG.$

The following  lemma is important for our study of random
generation. We will need to show that the homology of certain
cycles in a  Schreier graph of a $p$-group $K$ with 3 variables
$S=\{g_1,g_2,g_3\}$ is rich enough. We first need to make some
definitions.

\begin{definition}\lb{kappa}
For  $v\in X$, let $\kappa(v)$ denote the smallest positive
power of $p$ so that there exists a word $w$ in $\{g_1, g_2,
g_*= g_3^{\kappa}\}$ so that $v^{\overline w}=v$ and the
exponent sum of $g_*$ in $w$ is not divisible by $p$.  Let
$v_X$ denote a vertex for which $\kappa$ is minimal, let $w_X$
denote the corresponding word, and let $\kappa(X)=\kappa(v_X)$.
\end{definition}

Assume that
 \bel{ass0} \mbox{
 $g_*=g_3^{\kappa(X)}$ fixes no point in $X$}.
 \end{equation}
For a set of words $W$, and an element $f$ of the free group
$\mathcal{F}_W$ indexed by $W$, let $w_f$ denote the
concatenation of words $w\in W$ according to $f$.

\begin{lemma}[Words with special homology]\lb{wsh}  Assume
(\ref{ass0}). There exists a set $W$ of $|v_X^{\langle
g_1,g_2\rangle}|+1$ words in $S$ fixing $v_X$, so that for each
fixed $f\in \mathcal F_{W}$ and fixed $u\in X$ the 1-chains
$$ \tau(w_Xw_f)(\circ
u), \quad \left\{\tau w(v_X)\, :\, w\in W\cup \{w_X\}\right\}$$
are linearly independent.
\end{lemma}

\begin{proof}
Let $W$ be the generating set for the fundamental group of the
Schreier graph $\GG_{12}=\GG(K,v_X^{\langle
g_1,g_2\rangle},\{g_1,g_2\})$ defined in the proof of Lemma
\ref{pi1}, i.e. pick a spanning tree rooted at $v_X$, and for
each edge $(v_w,g_{i_w})$ outside the tree consider the path
$w(v_X)$ that moves from $v_X$ to $v_w$ in the tree, then along
the edge, and finally back to $v_X$ in the tree.

We need to check that if
$$
a_u \tau(w_Xw_f)(\circ u)+a_X\tau(w_X(v_X)) + \sum_{w\in W} a_w
\tau(w(v_X)) =0 \qquad (\mod p)
$$
then all the coefficients $a_i$ equal $0\ (\mod p)$. Indeed, by
Lemma \ref{linearind} below, the first two terms are linearly
independent when restricted to $g_3$ edges, and the terms in
the last sum vanish on $g_3$-edges, thus we have $a_+=a_u=0$.
It is standard that $\tau W$ are independent: they form a basis
for the homology group $\HH \GG_{12}$. Specifically, the edge
$(v_w,g_{i_w})$ defined in the previous paragraph appears in
exactly one of $\tau(v_X,w)$, so $a_w\equiv 0$.
\end{proof}

We used the following lemma in the proof above. Let $\tau_3
\pi$ denote the 1-chain of $\pi$ restricted to $g_3$-edges of
$\GG$. For a $\{g_1,g_2,g_*\}$-path $\pi$, let $\tau_* \pi$
denote the 1-chain in $\GG(X,\{g_1,g_2,g_*\})$ restricted to
$g_*$-edges.

If $w$ is word in which the exponent sum of $g_3$ is not
divisible by $p$, then it follows from Lemma \ref{homo} that
$\{\tau_3 w(u): u\in X\}$ are linearly independent.  The ideal
situation is when there is such a word satisfying $v^{\overline
w}=v$. But it can happen that there is no such word, so we need
the following lemma.

\newcommand{\vnull}{v}
\newcommand{\wcsillag}{w}
\begin{lemma}[Independence for added noise]\lb{linearind}
Assume (\ref{ass0}), and let $u\in X\setminus\{v_X\}$. The
vectors $\tau_3(w_X(v_X)), \; \tau_3(w_X(\circ u))$ are
linearly independent.
\end{lemma}

\begin{proof}
Let $\kappa=\kappa(X)$, $v=v_X$, $\wcsillag=w_X$, and let
$g_*=g_3^{\kappa}$. Consider the partition of $X$ into $\{
g_1,g_2,g_*\}$-orbits $\QQ$. For $x\in X$, let $Q_x\in \QQ$
denote the orbit of $x$.

Let $Q \in \QQ$. We claim that the sets $Q,\, Q^{g_3},\,
\ldots, Q^{g_3^{\kappa-1}}$ are disjoint.  Suppose the
contrary, then for $x,z\in Q$, we have $x^{g_3^i}=z^{g_3^j}$
with $i<j$. Then $x=z^{g_3^{j-i}}$. Let $w_{xz}$ be the word in
$\{g_1, g_2, g_*\}$ so that $x^{\overline w_{xz}}=z$. Then the
evaluation of the word $w_{xz}g_3^{j-i}$ fixes $z$ and its
$g_3$-exponent sum is not divisible by $\kappa$. But by
definition $\kappa$ is minimal at $v_X$, a contradiction.

Consider the homomorphism $\varphi$ of 1-chains defined by
\begin{eqnarray*}\varphi\;:\;\tau
\GG(Q,\{g_*\})&\longrightarrow &\tau \GG(X,\{g_3\}) \\
\tau_* g_*(x)&\mapsto& \tau_3 g_*(x)
\end{eqnarray*}
which maps the the 1-chain of a $g_*$-edge to the 1-chain of
its expansion, a path of length $\kappa$. By the previous
paragraph, for different $x,z\in Q$ the paths $g_*(x),\,
g_*(z)$ in $\GG(X,\{g_3\})$ have disjoint edges. It follows
that $\varphi$ is injective.

Since the exponent sum of $g_*$ in $\wcsillag$ is not divisible
by $p$, we can apply Lemma \ref{homo} to the Schreier graph
$\GG(Q,\{g_1,g_2,g_*\})$ and  $S'=\{g_*\}$. We get that
$\{\tau_*(\wcsillag(x)):\, x\in Q\}$ are linearly independent.
Since $\varphi$ is an injective homomorphism,
$\{\tau_3(\wcsillag(x)):\,x\in Q\}$ are linearly independent.

Let $x\in X$ and let $q$ denote the length of the cycle of $x$
for the action of $\overline \wcsillag$. Then
 \bel{nemnulla}
  \tau_3(\wcsillag(\circ x))=\sum_{\ell=0}^{q-1}
  \tau_3(\wcsillag(z^{{\overline \wcsillag}_*^\ell})) \not=0
 \end{equation}
since the terms are linearly independent. This implies that the
vectors $\tau_3(\wcsillag(\vnull)), \; \tau_3(\wcsillag(\circ
u))$ are nonzero. It also implies that if
 $u\in Q_{\vnull}$, then the claim of the lemma holds.

Now assume that $u\notin Q_{\vnull}$. Recall the definition of
the boundary operator $\partial$. It is a linear map from the
1-chains of a graph $\GG$ to 0-chains i.e. abstract linear
combinations of vertices of $\GG$. It is defined by a relation
$(x,z)\mapsto z-x$ for each edge $(x,z)$.

For $x\in X$,  $\tau_3 w(\circ x)$ is a linear combination of
1-chains of paths  with length $\kappa$ that start and end
$Q_{x}$. Thus $\partial \tau_3 w(\circ x)$ is supported on
$Q_x$. In particular, $\partial \tau_3 w(v)$, and $\partial
\tau_3 w(\circ u)$ are supported on $Q_v$, $Q_u$, respectively.
It therefore suffices to show that $\partial \tau_3
w(v)\not=0$.

We first claim that $\partial \tau_* w(\vnull)\not=0$, in other
words $\tau_* w(\vnull)$ is not a linear combination of
1-chains of $g_*$-cycles. Indeed, all such cycles are of
$p$-power length, and by assumption \eqref{ass0}, there is no
cycle of length $1$. But  $\tau_* w(\vnull)$ cannot be a linear
combination of longer cycles, since the exponent sum of $g_*$
in $\wcsillag$ is not divisible by $p$.

Since $\varphi$ is injective, it follows that $\partial \tau_3
w(\vnull)\not=0$, completing the proof.
\end{proof}

\section{Slices of random subgroups}\lb{srs}

The goal of this section is to show that slices of random
subgroups of $\Gamma(p)$ are large; in the next section we will
show that these large slices generate a large subgroup.

The setup is as follows. Let $K$ be a finite $p$-group acting
on a set $X$.

Let $H=\Gamma_n(p)$. The group $H\wr K$ acts on the disjoint
union of $|X|$ copies of $T_n$. For a subgroup $G_*\subseteq
(H\wr K)$, let $\sigma G_*\subseteq G_*$ denote the {\bf
boundary slice}, that is the pointwise stabilizer of the
vertices on level $n-1$ of the trees. Note that $\sigma G_*$ is
a vector space.

We first consider the case $K=1$. The following lemma shows
that it is possible to generate large slices by two elements,
even if one is required to be a high power.

\begin{lemma}\lb{polihamu}
For  $0\le k<n$ there exists elements $s,g\in \Gamma_{n}(p)$ so
that $$\dim (\sigma \langle s,g^{p^k} \rangle) =
p^{n-1}-p^k+1.$$
\end{lemma}

The following proof uses a technique introduced in \cite{av02},
but is self-contained.

\begin{proof}
Each element $v\in \partial T_{n-1}$ is represented as a
sequence $(a_1, \ldots, a_{n-1})$ of elements of $\FF_p$. We
assign to each such element a number $N(v)= a_1+a_2p+ \ldots+
a_np^{n-2}$. Recall that there exists an element, the adding
machine  $s\in \Gamma(p)$, that acts as a full cycle on
$\partial T_{n-1}$ so that $N(v^s)=N(v)+1\ (\mod p^{n-1})$.

Each $g\in \sigma \Gamma_{n}$ we identify with the polynomial
in $F_n:=\mathbb F_p[x]/(x^{p^{n-1}}-1)$ given by
$$
\sum_{v\in \partial T_{n-1}} g(v)x^{N(v)}
$$
this is an isomorphism of vector spaces, and by the above, the
action of $s$ on $\sigma\Gamma_{n}$ corresponds to
multiplication by $x$ in $F_n$. Define the {\bf weight}
$\nu(f)$ of an element $f\in F_n$ as the highest $k$ so that
$y^k=(x-1)^k$ divides $f$. Then it is clear that elements of
$F_n$ of different weight are linearly independent. Also
$\nu(xf-f)=\nu(yf)=\nu(f)+1$, unless $\nu(f)=p^{n-1}-1$ (since
$y^{p^{n-1}}=x^{p^{n-1}}-1=0$). This implies that if a
$s$-invariant subspace $M$ of $\sigma\Gamma_{n}$ contains an
element of weight $\nu$, then it also contains an element of
weight $\nu+1$, then one of $\nu+2$, up to weight $p^{n-1}-1$.
Hence $\dim M \ge p^{n-1}-\nu$.

Consider the element $g\in \Gamma_{n}$ with the following
properties. The action of $g$ agrees with $s$ on $\partial
T_k$. Also, for vertices $v$ at levels greater than $k$,
$g(v)=0$, except that for the vertex $v$ at level $n-1$ with
$N(v)=0$ we have $g(v)=1$.

Let $q=p^k$. One easily checks that $g^{q}\in
\sigma\Gamma_{n}$. Moreover, for $v\in \partial T_{n-1}$ we
have that $g^{q}(v)$ equals $1$ if $N(v)\le q$ and $0$
otherwise. Thus $g^{q}$ corresponds to the polynomial
$$f=1+x+ \ldots +x^{q-1}=y^{q-1}$$
of weight $q-1$ (the last equality can be checked via binomial
expansion). It follows that $g$ has weight $q-1$, and therefore
the $s$-invariant subspace containing it has elements of all
weights in $q-1\ldots p^{n-1}-1$. The claim of the lemma
follows.
\end{proof}

Now let $K=\langle \tilde g_1,\tilde g_2,\tilde g_3\rangle$.
Define $\kappa$, $v_X$, $w_X$ as in Definition \ref{kappa}
(with tildes removed), and assume $\ref{ass0}$. Consider the
subgroup $G\subseteq H\wr K$ generated by uniform random
liftings $g_i$ of the generators $\tilde g_i$. Let $\eps$
denote the probability that $|v^{\langle g_1,g_2\rangle}|+1$
random elements do not generate $H$.

\begin{lemma}[Controlled randomness in wreath products]
\lb{egylegeny}\ There exist a word $w$ in the $g_i$ whose
composition depends on the $g_i$ so that that $v_X^{\overline
w}=v_X$, $\pr({\overline w}(v_X)=h)\ge 1-\eps$, and for $u\in
X\setminus \{v_X\}$ the distribution of $\overline w(\circ u)$
is uniform on $H$.
\end{lemma}

\begin{proofof}{Lemma \ref{egylegeny}}
With the notation and the results of Lemma \ref{wsh}, consider
the group generated by $\{\overline w(v_X): w\in W\}$, with
probability $1-\eps$ this is the full group $H$. Then for some
concatenation  $w_f$ of elements of $W$ we have $\overline
w_+\overline w_f(v_X)=h$ (otherwise we pick $w_f$ arbitrarily).

The linear independence statement of Lemma \ref{wsh} translates
to probabilistic independence by Lemma \ref{lip}. Therefore,
for each fixed $f$, $u\not=v_X\in X$, given $\{\overline
w(v_X): w\in W\cup \{w_V\}\}$, $\overline w_+\overline w_f(u)$
has uniform distribution on $H$. This is also true for our
random choice of $f$ as it only depends on the values
$\{\overline w(v): w\in W_+\}$ (see Remark \ref{fuggetlen}).
\end{proofof}

\begin{remark}\lb{fuggetlen} In the last paragraph of the proof we used the
following simple fact. Let $X, \{Y_i,\, i\in I\}$ be random
variables so that for each $i\in I$, the two variables $X,\,
Y_i$ are independent and $Y_i$ has distribution $\mu$. Let
$f(x)$ be an $I$-valued deterministic function. Then
$X,Y_{f(X)}$ are also independent and $Y_{f(X)}$ has
distribution $\mu$.
\end{remark}
The {\bf rigid vertex stabilizer}  $\RR_{v}G\subseteq G$ of a
vertex $v$ is the pointwise stabilizer of the vertices not
descendant to $v$.


\begin{lemma}[Slices for a descendant tree of a vertex]\lb{becsiszelet4}
There exists $c_i>0$ depending on $p$
only so that if $\ell:=|v^{\langle g_1,g_2\rangle}|<p^{c_0n}$
then the event
 \bel{fatslice4}
 \dim \sigma \RR_{v_X}G \ge ( 1-p^{-c_1n+1}) \dim  \sigma \Gamma_n
 \end{equation}
has probability at least $1-e^{-c_2n}\ell-e^{c_3(n-\ell)}$.
\end{lemma}

\begin{proof} Let $v=v_X$, let
$G_v\subseteq G$ denote the stabilizer of $v$, let $G_{v+}=\{h
\in G_{v}: h(v)=1\}$, and let $\pi$ denote the quotient map
$G_v \rightarrow G_v / G_{v+} \equiv H$. Fix $c_1>0$ and let
$n_1=\lfloor c_1 n \rfloor$.

Consider the two elements $g,\,s$ used in Lemma \ref{polihamu}
for generating a large slice. Using Lemma \ref{egylegeny}
twice, with high probability we generate elements $g'\,s'\in
G_v$, so that $\pi g'=g$, $\pi s=g$. If $h'=g^{p^{n_1}}$ then
$h=\pi h'\in \sigma \Gamma_{n}$, and the smallest $s$-invariant
subspace $M$  of $\sigma \Gamma_n$ containing $h$ has
\bel{mdim} \dim M \ge  p^{n-1}-p^{n_1}+1 \end{equation}

We will show that with high probability, $h\in \RR_{v}G$, and
thus so are its $s$-conjugates. Therefore the subspace $M$ has
an isomorphic pre-image $M_v\subseteq G\cap {\RR_{v}G}$. The
inequality \eqref{mdim} implies the claim \eqref{fatslice4}.

Indeed, $h(u)$ for $u\in X\setminus\{v\}$ is the $p^{n_1-r}$-th
power of a uniform random element, where $p^r$ is the length of
the orbit of $u$ under $h$. By Theorem \ref{brw} if $c_5<1$ is
large enough, then there exist $c_2>0$ so that for all $n$ and
a uniformly chosen random $h_0\in H$ has $ \pr(h_0^{p^{\lfloor
c_5n \rfloor}} \not= 1) \le  e^{-c_2 n}$. With the right choice
of $c_0, c_1$, we get that
$$\pr\left(h(u)\not=1\mbox{\; for some\;
}u\in X\setminus\{v\}\right) \le p^m e^{-c_2n}.$$

The probability that Lemma \ref{egylegeny} does not give the
right $h$ is bounded above by the chance $q$ that $p^m+1$
random elements do not generate $\Gamma_n$ (or, equivalently,
its Frattini quotient $\mathbb \FF_p^n$). It is known that $q
\le e^{c_3(n-p^m)}$ where $c_3$ is an absolute constant. Since
we use Lemma \ref{egylegeny} twice, we can increase $c_3$ to
include a factor of $2$.
\end{proof}

Let $r$ denote the number of orbits of $K$ on $X$, and let
$\ell$ denote the size of the shortest orbit.

\begin{lemma}[Slices of random subgroups are large]\lb{becsiszelet3}
There exists $c_0,c_1,c_2,c_3>0$ depending on $p$ only so that
if $\ell<p^{c_0n}$ then the event
 \bel{fatslice}
 \dim \sigma G \ge ( 1-p^{-c_1n+1})\dim  \sigma (H \wr K)
 \end{equation}
has probability at least $1-r(e^{-c_2n}\ell-e^{c_3(n-\ell)})$.
\end{lemma}

\begin{proof}
For a $K$-orbit $R$ of $X$, let $\kappa(X)=\kappa(v_R)$ as in
Definition \ref{kappa}. Let $R_1,\ldots, R_r$ denote the
$K$-orbits of $X$, listed so that $\kappa(R_\ell)$ is
increasing in $\ell$. Let $X_\ell=R_\ell \cup\ldots \cup R_r$,
and let $D_\ell$ denote the subgroup of $D=H \wr K$ which fixes
all descendants of vertices in $X_\ell$, and let $G_\ell=D_\ell
\cap G$. Then we have
$$
1=D_1 \vartriangleright \ldots \vartriangleright D_{r+1}=D,
$$
and similarly
$$
1=G_1 \vartriangleright \ldots \vartriangleright G_{r+1}=G.
$$
Let $A=\sigma(\Gamma_n\wr K)$ be the
subgroup fixing all vertices at level $n-1$, a vector space. Then we
have
 \bel{osszeg}
  \dim (\sigma G) = \sum_{\ell=1}^r \dim \left(\sigma G_{\ell+1} \,/\, \sigma
  G_\ell\right).
 \end{equation}
Note that the subgroup $G/G_{\ell}\subseteq D/D_{\ell}$ acts
naturally on $X_\ell$, moreover, the $g_i/D_{\ell}$ are uniform
random liftings of the $\tilde g_i/(K\cap D_{\ell})$ to
$D/D_{\ell}$.

We can apply Lemma \ref{becsiszelet4} to this action with
$X=X_\ell$, $v_X=v_{R_\ell}$, since $v_{R_\ell}$ minimizes
$\kappa$ over $v_{X_{\ell}}$. We get a large boundary slice of
the subgroup $ \RR_v(D/D_{\ell})=\RR_v(D_{\ell+1}/D_\ell)$.
Conjugation by $G_{\ell+1}/G_{\ell}$ gives the $|R_\ell|$-fold
direct product of this large slice (one copy for each vertex
$u\in R_\ell$). Therefore with high probability
$$
\dim \left(\sigma G_{\ell+1} \,/\, \sigma G_{\ell}\right) \ge
|R_\ell| \; ( 1-p^{-c_1n+1})\,\dim\sigma \Gamma_{n}
$$
since Lemma \ref{becsiszelet4} is used $r$ times, the
probability of success is bounded below by the stated amount.
The claim of the Lemma now follows from (\ref{osszeg}) and the
straightforward identity $\dim \sigma(\Gamma_n \wr K)=|X| \dim
\sigma \Gamma_n$.
\end{proof}

\section{Dimension of random subgroups}\lb{drs}


The notion of Hausdorff dimension was introduced by
\cite{abercrombie94} to the setting of compact groups. Let $n$
be the greatest number so that $g\in \Gamma^{(n)}$. Define the
norm and distance
\bel{metric} \|g\|=|\Gamma_n|^{-1},\quad
\mbox{dist}(g,h):=\|gh^{-1}\|. \end{equation}  This turns
$\Gamma(H)$ into a compact metric space, and Hausdorff
dimension is defined in the usual way. Translated to our
setting the result of \cite{bs97} says that for closed
subgroups $G \subseteq \Gamma(H)$ Hausdorff dimension agrees
with lower Minkowski (box) dimension. Considering what balls
are in $\Gamma(H)$, this is easily checked to coincide with the
lim inf of the {\bf density sequence} $\gamma_\ell$ defined in
the introduction \eqref{density}. For technical purposes, we
will now introduce $\gamma_\ell$ in a slightly more general
setting, for the wreath product of $\Gamma(p)$ and a finite
$p$-group.

Let $K_0$ be a $p$-group acting on a set $X_0$. Let
$\Gamma_{K_0} = \Gamma(p) \wr K_0$ acting on $|{X_0}|$ copies
of the infinite $p$-ary tree $T_p$. Let $\Gamma_+^{(\ell)}$
denote the subgroup that stabilizes all vertices at level
$\ell$ of all trees. Given a subgroup $G\subseteq \Gamma_+$,
define the density sequence
$$
\gamma_\ell(G)=\frac{\log |G / \Gamma_+^{(\ell)}|}{\log
|\Gamma_+/ \Gamma_+^{(\ell)}|}.
$$
The density sequence $\gamma_\ell$ defined in \eqref{density}
covers the case when ${K_0}$ is the trivial group acting on a
singe-element set. In general, we have
$$\dimh(G)=\liminf_{\ell \rightarrow \infty} \gamma_\ell(G).$$

\begin{proposition}\lb{1d}
  Let $G$ be the the subgroup generated by independent uniform
random liftings of three fixed elements of ${K_0}$ to
$\Gamma_+=\Gamma(p)\wr K_0$. With probability 1,  the density
sequence $\gamma_\ell(G)$ satisfies
$$
\gamma_\ell(G) > 1- e^{-c\ell}
$$
for $c=c(p)>0$ fixed and all sufficiently large $\ell$. In
particular, $\dimh(G)=1$.
\end{proposition}

The most important case is when ${K_0}$ acts trivially on the
single-element set ${X_0}$.

\begin{theorem}\lb{3random} The subgroup of
$\Gamma(p)$ generated by three random elements has
1-dimensional closure with probability 1.
\end{theorem}

These random subgroups are the first known examples of finitely
generated subgroups of $\Gamma(p)$ with full dimensional
closure. We propose

\begin{conjecture}\lb{tworandom}
 Theorem \ref{3random} holds even for two random
elements.
\end{conjecture}

Another immediate corollary of  Proposition \ref{1d} is
\newtheorem*{dixon.thm}{Theorem \ref{dixon}}
\begin{dixon.thm}
\dixonthm
\end{dixon.thm}
This theorem is the $p$-group analogy of the famous result of
\cite{dixon69} saying that two random elements of $\Sym(n)$
generate $\Sym(n)$ or the alternating group $\Alt(n)$ with
probability tending to 1. In fact, for finite simple groups
$G$, the probability that two elements generate $G$ tends to
$1$ as $|G|\rightarrow \infty$, see \cite{shalev99}. The
methods used here are fundamentally different from the usual
subgroup counting technique.

Another interesting consequence of the Proposition \ref{1d} is
the following strengthening of Theorem \ref{tetszdim}.

\begin{theorem}\lb{mindend} For each $d\in [0,1]$ there exists a
topologically finitely generated free pro-$p$ subgroup of
$\Gamma(p)$ of Hausdorff dimension $d$.
\end{theorem}

This result leads to the solution of a problem of
\cite{MR2001d:20026}, who asked whether a topologically
finitely generated pro-$p$ group can contain topologically
finitely generated subgroups of irrational Hausdorff dimension.
Indeed, we may take three elements generating a
full-dimensional subgroup $G_1\subseteq \Gamma(p)$, and three
more elements that generate a $d$-dimensional subgroup
$G_d\subseteq \Gamma(p)$. It is straightforward to check that
the dimension of $G_d$ in $\langle G_1, G_d\rangle$ with
respect to the filtration inherited from $\Gamma(p)$ is $d$.

\newcommand{\Gammaprime}{{\Gamma_{T'}}}

Let $T'$ be a subtree of $T$ which has the same root as $T$ and
each vertex has $0$ or $p$ offspring. Let $\Gammaprime\subseteq
\Gamma(p)$ be the pointwise stabilizer of $T'$. The measure
$\mu(T')$ of $T'$ can be defined as $\lim |T'_n|/|T_n|$, where
$T'_n$ means the vertices at level $n$. The sequence is
non-increasing, so the limit always exists. It is easy to check
that
$$
\dimh(\Gammaprime) = 1- \mu(T').
$$
Let $V_{0}$ denote the set of vertices of $T'$ that have no
children. Since $\mu(T')$ can take any value in $[0,1)$ (even
if we assume that $V_0$ is infinite), $\dimh\ \Gammaprime$ can
take any value in $(0,1]$. In the following we show that random
subgroups  $G$ of $\Gammaprime$ are free and $\dimh(G)=
\dimh(\Gammaprime)$ with probability 1. This proves Theorem
\ref{mindend} for $d>0$; the easy case $d=0$ is left to the
reader.

\begin{proposition}\lb{spectrumpro} The subgroup
$G$ topologically generated by $k\ge 3$ random elements of
$\Gammaprime$ satisfies
$$ \dimh (G) = \dimh(\Gammaprime) = 1-\mu(T')
$$
with probability 1. If $V_0$ is infinite, then $G$ is a free
pro-$p$ group with probability 1.
\end{proposition}

\begin{proofof}{Proposition \ref{spectrumpro}} Let $V_{0+}$
denote the vertices $V_{0}$ and their descendants
in $T$. Clearly, for $g\in \Gammaprime$ and $v\in V\setminus
V_{0+}$ we have $g(v)=0$. Moreover, Haar measure on
$\Gammaprime$ agrees with picking $g(v)$ uniformly,
independently at random for each $v\in V_{0+}$.

Let $g_i$ denote the random elements, let $\eps>0$, and let
$\ell$ be so that $|\partial_\ell T'|/|\partial_\ell
T|<\mu(T')+\eps$. Condition on the action ${K_0}$ of $g_i$ on
${X_0}=\partial T_\ell\setminus
\partial T'_\ell$. By the first paragraph, given this information,
the $g_i$ are uniform random liftings. Proposition \ref{1d}
applied to this action, with dimension appropriately rescaled,
implies
$$
\dimh G \ge |\partial _\ell T \setminus \partial_\ell
T'|/|\partial_\ell T|>1-\mu(T')-\eps.
$$
With probability 1, This is then true for all rational $\eps$,
so $\dimh(G) \ge \dimh(\Gammaprime)$, and the other inequality
is trivial.

For the vertices $v\in V_0$ the action of $G$ on the descendant
subtrees of $v$ are independent, generated by $k$ independent
random elements. If $V_0$ is infinite, then each finite
$p$-group that is generated by $k$ elements appears with
probability 1 as a quotient group of the action on one of the
descendant subtrees. Thus each finite $p$-group generated by
$k$ random elements appears as a quotient group of $G$ with
probability 1, which implies that $G$ is a free pro-$p$ group
with probability 1.
\end{proofof}

We now proceed to prove the key proposition. The proof uses
almost all theorems proved so far, in particular, it depends
crucially on the results of Sections \ref{sgh} and \ref{srs}.

\begin{proofof}{Proposition \ref{1d}}
For $\ell \ge 1$ let $m=m(\ell):=\lfloor (\log \ell)^2\rfloor$,
let $n:=\ell-m$, and condition on the values of
$g_i':=g_i/\Gamma_+^{(m)}\in \Gamma_m \wr {K_0}$. Given this
information, the $g_i/\Gamma^{(\ell)}$ are uniform random
liftings of $g_i'$ to $\Gamma_n \wr  (\Gamma_m \wr
{K_0})=\Gamma_\ell \wr {K_0}$. Let $X_m$ denote the set of
vertices at level $m$ in the trees. Let $b\ge 1$, and assume
that the following event $A_{b,m}$ holds. Let $G_{12}=\langle
g_1, g_2 \rangle$. For the action of $G$ on $X_m$:
 {\list{}{\itemsep-.5ex \topsep .5ex}
 \item[(i)] the number of orbits of $(G_{12},X_m)$ is at most
 $b$,
 \item[(ii)] the shortest orbit of $(G_{12},X_m)$ has length at
 least $p^m/b$,
 \item[(iii)] if $p^m>b$, then $g^{\kappa(X_m)}_3$ has no fixed
 points.
 \endlist}

We now apply Lemma \ref{becsiszelet3} with $X=X_m$, $K=\Gamma_m
\wr {K_0}$ here. Let $B_\ell$ be the event that
$$
s_\ell:=\dim\left((G\cap
\Gamma_+^{(\ell-1)})\,/\,\Gamma_+^{(\ell)}\right)
> p^{\ell-1}(1-p^{-c_6 \ell}),
$$
call this event $B_\ell$. Lemma \ref{becsiszelet3} implies
$\pr(B_\ell^c|A_{b,m})\le e^{-c_4\ell+c_5}$ for some constants
$c_4,\,c_5>0$ depending on $b$, $p$.

By Proposition \ref{kopaszfa}, for all large enough $m$ the
children of any $G_{12}$-orbit of $X_m$ form a $G_{12}$-orbit
of $X_{m+1}$. Therefore assumptions (i)-(ii) hold for some
random $b$ and all $m$ with probability 1. Also, we get that if
$w$ is a word in the generators fixing $v\in X_m$, then for
some $\{g_1,g_2\}$-word $w'$ all children of $v$ are fixed by
$ww'$. Therefore $\kappa(X_{m+1}) \le \kappa(X_{m-1})$, so
$\kappa(X_m)<b''$ for some random $b''$ and all $m$.

Every random word in $\Gamma(p)$ fixes only finitely many fixed
vertices of $T$ (Proposition \ref{wordfix}), so (iii) holds for
some random $b$ and all $m$ with probability 1.
 Thus if we define
$$
A_b=\bigcap_{m \ge 1} A_{b,m}
$$
then  $A_b$ is increasing in $b$ and $\pr(A_b) \rightarrow 1$.

Now we have
$$
\pr(B_\ell \cap A_b)\le \pr(B_\ell \cap A_{n,b}) \le
\pr(B_\ell|A_{n,b})
$$
Since the latter is summable, we get that for fixed $b$,
$B_\ell\cap A_b$ happens for at most finitely many $\ell$ with
probability 1.  Taking the increasing union as $b\rightarrow
\infty$ we get that $B_\ell$ happens for at most finitely many
$\ell$ with probability 1.  Then
\begin{eqnarray*} \log_p
|G/\Gamma_+^{(\ell)}| &=& \sum_{k=1}^{\ell} s_k
> \sum_{k=\ell_0}^{\ell}  |{X_0}|\, p^{\ell-1}(1-p^{-c_6 \ell})
\\&>& |{X_0}|\, \frac{p^\ell-1}{p-1} (1-e^{-cn}) = \log_p
|\Gamma_+/\Gamma_+^{(\ell)}|(1-e^{-cn})
 \end{eqnarray*}
 where the last
inequality holds for large enough $\ell$ if $c<c_6\log p$. The
claim of the proposition follows.
\end{proofof}

We conclude this section with two conjectures of increasing
strength.

\begin{conjecture}\lb{bigfg} Let $G$ be a closed subgroup
of $\Gamma(p)$. Then $G$ contains a topologically finitely
generated subgroup of the same Hausdorff dimension.
\end{conjecture}

It may be that such subgroups are abundant. An affirmative
answer to the following question would imply Conjectures
\ref{bigfg} and \ref{tworandom}.

\begin{question} Let $G$ be a closed subgroup of $\Gamma(p)$. Let
$g_1, g_2$ be random elements chosen according to Haar measure
on $G$. Is it always true that  $\dimh(\langle g_1,
g_2\rangle)=\dimh(G)$ with probability 1?
\end{question}

\section{Small subgroups}
\lb{small}

In this section we state general results on small subgroups of
$\Gamma (p)$.

First we give a description on the Abelian subgroups of $\Gamma
(p)$ using orbit trees. A rooted tree is a $1$-$p$ tree if all
the vertices have $1$ or $p$ children. We call a vertex {\bf
solo}, if it has $1$ child. For an
arbitrary rooted tree $R$ let $S(R)$ denote the number of solo vertices of $%
R $. We say that a group $A\subseteq \Gamma _{n}(p)$ {\bf
belongs to} $R$ if $R=T_{n}/A$.

\begin{lemma}
\lb{egyes}Let $R$ be a $1$-$p$ tree of depth $n$ and let
$A\subseteq \Gamma _{n}(p)$ an Abelian subgroup which belongs
to $R$. Then $\log _{p}\left| A\right| \leq S(R)$. Equality
holds if and only if $A$ is a maximal Abelian subgroup which
belongs to $R$.
\end{lemma}

\begin{proof}
We use induction on $n$. For $n=1$ the proposition is trivial.
Suppose it holds for $n-1$.

Let $B\subseteq A$ be the stabilizer of level 1 of $T_{n}$.
Then $B$ has index $1$ or $p$ in $A$. Let $\Lambda _{1},\ldots
,\Lambda _{p}$ denote the subtrees of $T_{n}$ starting from a
vertex of level $1$. Let $B_{i}$ denote
the restriction of the action of $B$ to $\Lambda _{i}$ and let $
\left| B_{i}\right| $.

If $A=B$ (i.e., if the root of $R$ is not solo), then $A$ maps
each $\Lambda _{i}$ onto itself so we have
\begin{equation*}
\left| A\right| \leq \prod_{i=1}^{p}\left| B_{i}\right|
=p^{\sum_{i=1}^{p}k_{i}}\text{.}
\end{equation*}
Using induction we get
\begin{equation*}
\sum_{i=1}^{p}k_{i}\leq \sum_{i=1}^{p}S(\Lambda
_{i}/B_{i})=S(R).
\end{equation*}

If $\left| A:B\right| =$ $p$ (i.e., the root is solo), let
$g\in A\setminus
B $ be an arbitrary element. Then $g^{p}\in B$, so the action of $g$ on $%
T_{n}/B$ has order $p$. This means that $g$ induces graph
isomorphisms between the trees $\Lambda _{i}/B_{i}$, so
$S(T_{n}/B)=pS(\Lambda _{1}/B_{1}) $ and $S(R)=S(\Lambda
_{1}/B_{1})+1$. Now let $b\in B$ be any
element which lies in the kernel of the restriction of $B$ to $\Lambda _{1}$%
, i.e., fixes $\Lambda _{1}$ pointwise. Since $g$ permutes the
$\Lambda _{i}$ transitively and commutes with $b$, it means
that $b$ fixes all the $\Lambda
_{i}$ pointwise, that is, $b$ is the identity. Thus $B$ is isomorphic to $%
B_{1}$. So using induction
\begin{equation*}
\left| A\right| =p\left| B_{1}\right| \leq p^{S(\Lambda
_{1}/B_{1})+1}=p^{S(R)}\text{.}
\end{equation*}

If the equality $\left| A\right| =p^{S(R)}$ holds, then $A$ is
trivially maximal with the property $R=T_{n}/A$. The other
direction follows by induction the same way as the inequality
above.
\end{proof}

As a corollary, one can show that Abelian subgroups of $\Gamma
(p)$ are zero-dimensional. Theorem \ref{solsum} gives a more
general result in this direction. Moreover, Theorem
\ref{perfekt} shows that Abelian groups are also
small-dimensional as sections: the derived subgroup of a closed
subgroup $G\subseteq \Gamma (p)$ has the same dimension as $G$.

For Theorem \ref{perfekt} we need a general asymptotic result
on rooted
trees. Let $d\geq 2$ be an integer. A rooted tree $T$ is defined to be a $d$%
-bounded tree if every vertex has at most $d$ children. It is a
$1$-$d$ tree
if every vertex has either $1$ or $d$ children. Let $U$ be a subset of a $d$%
-bounded tree $T$. We define the \emph{density sequence} of $U$
by
\begin{equation*}
\delta _{n}(U)=\frac{r_{n}(U)}{d^{n}}
\end{equation*}
where $r_{n}(U)$ is the number of vertices in $U$ of level $n$.
We define the {\bf density}
\begin{equation*}
\delta (U)=\lim_{n\rightarrow \infty }\delta _{n}(U)
\end{equation*}
if this limit exists. In particular, $\delta (T)$ always exists, since $%
\delta _{n}(T)$ is monotone decreasing.

The following straightforward lemma will be useful.

\begin{lemma}
\lb{falemma}Let $T$ be a $d$-bounded tree and let $U$ be the
set of vertices of degree less than $d$. Then
$\sum_{i=0}^{\infty }\delta
_{i}(U)\leq d(1-\delta (T))$. In particular, $\delta (U)=0$. If in addition $%
T$ is a $1$-$d$ tree then
\begin{equation*}
\sum_{i=0}^{\infty }\delta _{i}(U)=\frac{d}{d-1}(1-\delta
(T))\text{.}
\end{equation*}
\end{lemma}

\begin{proof}
For the first claim, we may assume that all the vertices have degree $d$ or $%
d-1$. Then $r_{i}(U)=dr_{i}(T)-r_{i+1}(T)$. From this, using
$r_{0}(T)=1$ we have
\begin{equation*}
\sum_{i=0}^{n}\delta _{i}(U)=\sum_{i=0}^{n}d\left( \frac{r_{i}(T)}{d^{i}}-%
\frac{r_{i+1}(T)}{d^{i+1}}\right) =d\left( 1-\frac{r_{n+1}(T)}{d^{n+1}}%
\right)
\end{equation*}
which gives $\sum_{i=1}^{\infty }\delta _{i}(U)=d(1-\delta
(T))$. If every vertex has degree $1$ or $d$ then we have
$r_{i}(U)=\left( dr_{i}(T)-r_{i+1}(T)\right) /(d-1)$ and the
last statement of the lemma follows.
\end{proof}

Now we are ready to prove the ``perfectness'' theorem, restated
here.

\newtheorem*{perfekt.thm}{Theorem \ref{perfekt}}

\begin{perfekt.thm}\perfektthm
\end{perfekt.thm}

\begin{proof}
Let $F=T(p)/G$ be the orbit tree of the action of $G$ on
$T(p)$. We treat $F$ and $T(p)$ as $1$-$p$ trees. Fix positive
integers $d$ and $k$.

We partition the vertices of $F$ into three classes as follows.
Let $A$ be the set of orbits which have exactly $p^{k}$ $k$-th
cousins in $F$ including themselves. Let $B$ be the set of
orbits of size at least $p^{d}$ and let $C$ be the set of
remaining vertices.

We claim that $\delta (B\cup C)=0$. To see this, let us define a new graph $%
F^{\prime }$ on the vertices on $F$ by putting an edge between
$u$ and $v$ if $u$ is the $k$-th grandson of $v$. Then
$F^{\prime }$ is the union of $k$
disjoint $p^{k}$-bounded trees. Now let $D$ be the set of $F^{\prime }$%
-parents of elements in $B\cup C$. Then $D$ is the set of vertices in $%
F^{\prime }$ which has less than $p^{k}$ children, so by Lemma
\ref{falemma}
$D$ has zero density in each of the components of $F^{\prime }$. But then $%
B\cup C$ has zero density in $F$.

For an arbitrary $p$-group $P$ acting on a finite set $X$ of
size $m$ let us define the {\bf commutator density} by
$c(P)=\log _{p}\left| P:P^{\prime }\right| $. We use the
following facts on $c(P)$.

\begin{itemize}
\item[(a)]  Comparing $P$ to the Sylow $p$-subgroup we have $c(P) \leq
\log_p \left| \mbox{Syl}_{p}(\mbox{Sym}(m))\right| \leq m$.

\item[(b)]   If $P\;$is transitive on $X$ then $c(P)\leq Km/\sqrt{\log m}$
where $K$ is an absolute constant. This is a result of
\cite{MR89i:20008}.

\item[(c)]  If $P$ is intransitive and $X=\bigcup_{i=1}^{j}X_{i}$ such that
each $X_{i}$ is $P$-invariant, then we have $c(P)\leq
\sum_{i=1}^{j}c(P_{i})$ where $P_{i}$ is the action of $P$ on
$X_{i}$.

\item[(d)]  If $P$ acts diagonally on the $X_{i}$, i.e., when elements of $P$
fixing any of the $X_{i}$ pointwise must fix the whole $X$, then $%
c(P)=c(P_{1})$ .
\end{itemize}

Now let $n$ be an arbitrary level and let $A_{n}$, $B_{n}$ and
$C_{n}$ be the set of vertices of $T$ which belong to $A$, $B$
and $C$ respectively. Let $G_{n}$ be the action of $G$ on the
$n$-th level of $T$ and let $G_{a}$, $G_{b}$ and $G_{c}$ be the
action of $G_{n}$ on $A_{n}$, $B_{n}$ and $C_{n}$.

The orbits lying in $A$ can be partitioned into sets of size
$p^{k}$ such that the action of $G_{a}$ on them is diagonal.
Then (d), (c) and (a) imply
\begin{equation*}
c(G_{a})\leq \left| A_{n}\right| /p^{k}\leq p^{n-k}.
\end{equation*}
Using (b) and (c) we have
\begin{equation*}
c(G_{b})\leq K\left| B_{n}\right| /\sqrt{d}\leq
Kp^{n}/\sqrt{d},
\end{equation*}
since all the orbits lying in $B$ have size at least $p^{d} $.
Lastly (a) implies
\begin{equation*}
c(G_{c})\leq \left| C_{n}\right| \leq \delta _{n}(C)p^{n+d-1},
\end{equation*}
since all the orbits in $C_{n}$ have size at most $p^{d-1}$.
Combining the above, (c) on $A_{n}$, $B_{n}$ and $C_{n}$
implies
\begin{equation*}
c(G_{n})p^{-n} \leq p^{-k}+K/\sqrt{d}+\delta _{n}(C)p^{d-1}
\end{equation*}
and letting $n\rightarrow \infty$ we get \bel{this}
\limsup_{n\rightarrow \infty} c(G_{n}) p^{-n} \leq p^{-k}+K/\sqrt{d}%
+\delta(C) p^{d-1}
\end{equation}
where $\delta (C)\leq \delta (B\cup C)=0$. The inequality
\eqref{this} holds for every positive $k$ and $d$, so
$c(G_{n})p^{-n}\rightarrow 0$, and the proof is complete.
\end{proof}

As an immediate corollary, solvable subgroups are
zero-dimensional. Using a Baire category argument we obtain the
following.

\begin{corollary}
\lb{sokfelold}Let $G\subseteq \Gamma (p)$ be a
positive-dimensional subgroup. Then $G$ cannot be abstractly
generated by countably many solvable subgroups.
\end{corollary}

\begin{proof}
Suppose that $H_{1},H_{2},\ldots \subseteq G$ are solvable
subgroups generating $G$ as an abstract group. Since closure
preserves solvability, we can assume that the $H_{i}$ are
closed. Then
\begin{equation*}
G=\bigcup_{k,n_{i}\in \mathbb{N}}H_{n_{1}}H_{n_{2}}\cdots
H_{n_{k}}\newline
\end{equation*}
is the countable union of the finite products formed from the
$H_{i}$. The sets $H_{n_{1}}H_{n_{2}}\cdots H_{n_{k}}$ are
closed, so by the Baire category argument one of them has to
contain an open set $U$. Then $U$ contains a coset of an open
subgroup in $G$, so $\dimh H_{n_{1}}H_{n_{2}}\cdots
H_{n_{k}}\geq \dimh U=\dimh G>0$.

By Theorem \ref{perfekt} $H_{n_{i}}$ are zero-dimensional,
moreover, $\lim \gamma _{n}(H_{n_{i}})=0$. It is easy to see
that
\begin{equation*}
\gamma _{n}(H_{n_{1}}H_{n_{2}}\cdots H_{n_{k}})\leq
\sum_{i=1}^{k}\gamma _{n}(H_{n_{i}})
\end{equation*}
for all $n$. This implies $\dimh H_{n_{1}}H_{n_{2}}\cdots
H_{n_{k}}=0$, a contradiction.
\end{proof}

Note that the full $\Gamma (p)$ can be easily generated
\emph{topologically} by an element and an Abelian subgroup.

Since solvable subgroups have dimension zero, we need to find a
more refined way to measure how large they are. This can be
done by \emph{summing} the density sequence rather than taking
its limit, as Theorem \ref{solsum} shows.

\newtheorem*{solsum.thm}{Theorem \ref{solsum}}
\begin{solsum.thm}
\solsumthm
\end{solsum.thm}

The example $\Gamma_n\subset \Gamma(p)$ shows that this is the
best possible bound up to the constant factor.

\begin{proof}
We can directly deduce this theorem from Lemma \ref{egyes} and
Lemma \ref {falemma} in the case when $G$ is Abelian.

Instead of $\gamma _{n}(G)$, it is more convenient to first
estimate a slightly different density. Let
\begin{equation*}
\overline{\gamma }_{n}(G)=(p-1)\frac{\log \left| G_{n}\right|
}{p^{n}}, \qquad s_n(G)=\sum_{i=0}^n \overline{\gamma }_i(G)
\text{. }
\end{equation*}

Let $m_{n,d}=\max \left\{ s_n(H)\mid H\subseteq \Gamma
_{n}(p)\text{ with solvable length }d\right\} $.

Let $\Lambda _{1},\ldots ,\Lambda _{p}$ denote the subtrees of
$\Gamma _{n}$ starting from level $1$. Let $K\subseteq H$
denote the stabilizer of level $1
$. Either $K=H$ or $\left| H:K\right| =p$. Let $K_{j}$ denote the action of $%
K$ on $\Lambda _{j}$. The groups $K_{j}$ are $d$-solvable so
$s_{n-1}(K_j) \leq m_{n-1,d}$ for all $j$.

If $H=K$ then we have \bel{rhs1} s_n(H)\leq \sum_{j=1}^{p} s_
{n-1}(K_{j})/p\leq m_{n-1,d}.
\end{equation}
If $\left| H:K\right| =p$ then let $A$ denote the restriction
of $K$ to the set $\Lambda _{2}\cup \ldots \cup \Lambda_{n}$ and let $%
N\subseteq K_{1}$ be the kernel of this restriction, i.e., the
set of elements which act trivially outside $\Lambda _{1}$. Let
$h\in H\backslash K$ be an element moving the first level and
let $D=\left[ N,h\right] $. Then $h$
moves all the vertices of the first level so the restriction of $D$ to $%
\Lambda _{1}$ is again $N$. Hence from $D\subseteq G^{\prime }$ we see that $%
N$ is $d-1$-solvable. Now $\left| H\right| =p\left| N\right|
\left| A\right| $ and $\left| A\right| \leq |K_2|\cdots |K_p|$
so we have
\begin{equation*}
\overline{\gamma }_{i}(H)\leq \frac{1}{p^{i}}+\frac{1}{p}\left( \overline{%
\gamma }_{i-1}(N)+\sum_{j=2}^{p}\overline{\gamma
}_{i-1}(K_{j})\right),
\end{equation*}
from which we get \bel{rhs2} s_n(H)\leq
\frac{p}{p-1}+\frac{1}{p}\left(
m_{n-1,d-1}+(p-1)m_{n-1,d}\right) \text{.}
\end{equation}

Summarizing, we have $m_{n,0}=0$ and a recursive upper bound on
$m_{n,d}$
given by the maximum of the right hand side of \eqref{rhs1} and \eqref{rhs2}%
. From this the bound $m_{n,d}\leq dp^{2}/(p-1)$ easily
follows. On the
other hand, trivially $\left| \gamma _{n}(G)-\overline{\gamma }%
_{n}(G)\right| \leq 1/p^{n}$ which yields
\begin{equation*}
\sum_{n=0}^{\infty }\gamma _{n}(G)\leq Cd
\end{equation*}
\newline
with $C=(p^{2}+p)/(p-1)$.
\end{proof}

A possible generalization of the result on the
zero-dimensionality of solvable groups would be that every
subgroup of $\Gamma (p)$ which satisfies a nontrivial identity
is zero-dimensional. Note that Conjecture \ref{propsub} would
imply this even for pro-$p$ identities. Another possible
direction arises from the observation that linear groups tend
to be zero-dimensional in $\Gamma (p)$.

\begin{conjecture}
\lb{linzero}Let $R$ be a commutative pro-$p$ ring and let
$G\subseteq
GL_{n}(R)$ be a linear group over $R$. Then for every embedding of $G$ into $%
\Gamma (p)$ the image of $G$ has zero dimension.
\end{conjecture}

It is easy to see that the conjecture holds for $R=\mathbb{Z}_{p}$; in
this case $G$ is $p$-adic analytic and so is the product of finitely many
procyclic subgroups (\cite{analyticpropgroups}). This question is also
related to Conjecture \ref {propsub} in the sense that by
\cite{MR2000d:20031} nonabelian free pro-$p$ groups are not linear over
local fields, so Conjecture \ref{propsub} implies Conjecture \ref{linzero}
in the case when $R$ is a local field.

\section{Large subgroups} \lb{large}

In this section we analyze $1$-dimensional subgroups of $\Gamma
(p)$. The main results are that spherically transitive
$1$-dimensional groups have normal spectra $\left\{ 0,1\right\}
$ and that they contain nonabelian free pro-$p$ subgroups as
well as dense free subgroups. Recall that the normal spectra of
a group $G$ is the set of possible Hausdorff dimensions of
normal subgroups of $G$.

One of the natural sources of interesting subgroups in $\Gamma
(p)$ are just infinite pro-$p$ groups. A group $G$ is just
infinite if every proper quotient of $G$ is finite; $G$ is
hereditarily just infinite if all subgroups of $G$ of finite
index are just infinite. Just infinite pro-$p$ groups are
regarded as the simple groups in the pro-$p$ category. Many
questions on pro-$p$ groups can be reduced to the just infinite
case, since (unlike in the profinite setting \cite{MR1894095})
every finitely generated pro-$p$ group possesses just infinite
pro-$p$ quotients.

Building on the work of \cite{MR43:338}, \cite{MR2002f:20044}
showed that just infinite pro-$p$ groups fall into the
following two classes:
\begin{itemize}
\item[(i)]  groups containing an open normal subgroup which is the direct power
of a hereditarily just infinite pro-$p$ group;
\item[(ii)]  pro-$p$ branch groups.
\end{itemize}
Groups from class (ii) have a natural action on
$T(p)$. In fact, one can define them as spherically transitive subgroups $%
G\subseteq \Gamma (p)$ such that the rigid level stabilizers of
$G$ have finite index in $G$. The first example for such a
group was the closure of the first Grigorchuk group, a subgroup
of $\Gamma (2)$ with many remarkable properties.
\cite{MR2001h:20046} calculated that it has Hausdorff dimension
$5/8$.

In \cite{MR2001b:20001}, Boston suggested a direct connection
between Grigorchuk's classes and Hausdorff dimension: he
conjectured that a just infinite pro-$p$ group is branch if and
only if it can be obtained as a positive dimensional subgroup
of $\Gamma (p)$. In this section we contrast known properties
of branch groups with new results on $1$-dimensional groups.
The first result concerns the normal subgroup structure of such
groups.

In general, the normal spectra of $1$-dimensional groups can
easily be the full interval $\left[ 0,1\right] $. The simplest
example for this is the stabilizer of an infinite ray in
$\Gamma (p)$; it is the full countable direct power of $\Gamma
(p)$ with weighted dimensions $1/p^{n}$. However, if we assume
that the group is spherically transitive, i.e., it acts
transitively on every level, we get a completely different
picture. It turns out that these groups, just as the free
pro-$p$ group (see \cite {MR2001c:16078}), are ``simple'' in
terms of Hausdorff dimension.

\begin{theorem}
Let $G\subseteq \Gamma (p)$ be a spherically transitive
$1$-dimensional closed subgroup. Then for every $1\neq
N\vartriangleleft G$ we have $\dimh N=1$.
\end{theorem}

We need some technical lemmas. The first one is an asymptotic
version of
Theorem \ref{egyszeru}. As the existence of large Abelian subgroups in $%
\Gamma _{n}(p)$ shows, we have to assume that the density is
close to $1$.

\begin{lemma}
\lb{aszperf}There exist a constant $C$ depending only on $p$
such that for every $G\subseteq \Gamma _{n}(p)$ with $\gamma
_{n}(G)\geq 1-\epsilon $ we have $\gamma _{n}(G^{\prime })\geq
1-\epsilon -\delta$ where $\delta =C \min
\{n,-\log_{p}\epsilon\}^{-1/2}. $
\end{lemma}

\begin{proof}
We consider the action of $G$ on level $n$ of $T_{n}$. Let
$r_{i}$ be the number of orbits of size $p^{i}$ $(0\leq i\leq
n)$, and let $r$ denote the total number of orbits. Then
\bel{rp} \sum_{i=0}^{n}r_ip^i=p^{n}
\end{equation}
and
\begin{equation*}
\log _{p}\left| G\right| \leq \sum_{i=0}^{n}r_{i}\log
_{p}\left| \Gamma _{i}(p)\right|
=\sum_{i=0}^{n}r_i\frac{p^{i}-1}{p-1}=\frac{p^n-r}{p-1}.
\end{equation*}
\newline
From this and the density assumption
\begin{equation*}
\log _{p}\left| G\right| \ge (1-\epsilon)\frac{p^n-1 }{p-1}
\end{equation*}
we get the inequality $r \leq \epsilon (p^{n}-1)+1$. Using the
theorem of \cite{MR89i:20008} with (b) and (c) from the proof
of Theorem \ref{perfekt} we see that
\begin{equation*}
\log _{p}\left| G:G^{\prime }\right| \leq K\sum_{i=0}^{n}r_{i}\frac{p^{i}}{%
\sqrt{i}}
\end{equation*}
$\newline $where $K$ is an absolute constant. Separating the
sum at an arbitrary $m\leq n$ and using \eqref{rp} we get
\begin{equation*}
\sum_{i=0}^{n}\frac{r_i p^{i}}{\sqrt{i}} \le \frac{rp^m}{\sqrt{m}} +\frac{p^n%
}{\sqrt{m+1}}.
\end{equation*}
\newline
Therefore with $m_0= \min \left\{ n,-\log _{p}\epsilon \right\}$ and $%
C=4K(p-1)$ we get
\begin{equation*}
\log _{p}\left| G:G^{\prime }\right| \leq K(p^n-1)\frac{3+\epsilon p^{m}}{%
\sqrt{m}}\leq C\log_p|\Gamma_n(p)|/\sqrt{m_0},
\end{equation*}
\newline
and this implies the stated bound on the density $\gamma
_{n}(G^{\prime })$.
\end{proof}

Consider a subgroup $G\subseteq \Gamma _{n}(p)$ where $n$ may
be infinite. We define the {\bf rigid stabilizer} $\RR_k G$ of
level $k$ as follows. Recall that the rigid vertex stabilizer
$\RR_v G$ of a vertex $v$ is the pointwise stabilizer of the
vertices not descendant to $v$. Then $\RR_{k}G$ is defined as
the group generated by the $\RR_{v}G$ of vertices $v$ at level
$k$.


We need a slight reformulation of Theorem 4. of
\cite{MR2002f:20044}. This is the key theorem to prove that the
first Grigorchuk group, or more generally, torsion branch
groups are just infinite.

\begin{lemma}
Let $G\subseteq \Gamma _{n}(p)$ be transitive and let $g\in G$
be an element which moves a vertex at level $k$. Then the
normal subgroup generated by $g$ in $G$ contains
$(\RR_{k+1}G)^{\prime }$, the derived subgroup of the rigid
stabilizer of level $k+1$.
\end{lemma}

The proof is exactly the same as in \cite{MR2002f:20044}.

We also need a lemma which measures how close subdirect
products of high density are to direct products in size.

\begin{lemma}
\lb{szubnagy}Let $G\subseteq \Gamma_{n}(p)^{d}$ satisfying
$|G|\geq |\Gamma_n|^{d-\epsilon}$. Then $G$ contains a direct
product $H=H_{1}\times \cdots \times H_{d}$ with $|H| \ge
 |\Gamma_n|^{ d(1-\epsilon )}$.
\end{lemma}

\begin{proof}
Let $i\leq d$ be a coordinate. Project $G$ to the coordinates not equal to $%
i $ and let $H_{i}$ be the kernel of this projection. Then the
image has size at most $|\Gamma_n|^{d-1}$ so $|H_{i}|\ge
|\Gamma_n|^{1-\epsilon}$. Since the full direct product
$H=H_{1}\times \cdots \times H_{d}$ lies in $G$, the statement
of the lemma follows.
\end{proof}

Now we are ready to prove Theorem \ref{egyszeru}, restated
here.

\newtheorem*{egyszeru.thm}{Theorem \ref{egyszeru}}
\begin{egyszeru.thm}
\egyszeruthm
\end{egyszeru.thm}

\begin{proof}
Let $N\vartriangleleft G$ be a nontrivial normal subgroup, and
fix $k$ so that $N$ moves a vertex at level $k-1$. For $n\ge k$
let $G_n=G\Gamma^{(n)}/ \Gamma^{(n)}$, let $N_n=N\Gamma^{(n)}/
\Gamma^{(n)}$. Then $N_n$ is a nontrivial normal subgroup of
$G_n \subseteq \Gamma_n$. Let $K\subseteq G_n$ denote the
stabilizer of level $k$. Then $| G_n:K| \leq \left| \Gamma
_{k}\right| $ so by the 1-dimensionality of $G$
\begin{equation*}
\log |K|=\log |\Gamma_n|\, (1-\epsilon_n)
\end{equation*}
where $\epsilon_n \rightarrow 0$ as $n\rightarrow \infty$. Now
$K$ acts on the subtrees starting at level $k$ so it is
naturally a subgroup of $\Gamma _{n-k}(p)^{p^{k}}$. Lemma
\ref{szubnagy} implies that $K$ contains a large direct
product; in particular, the rigid stabilizer $\RR_{k}=\RR_k G$
satisfies
\begin{equation*}
\log |\RR_k| \ge p^k \log |\Gamma_{n-k}| \,(1-p^k\epsilon_n)
\ge \log |\Gamma_n|\, (1-\epsilon^{\prime}_n)
\end{equation*}
with $\epsilon^{\prime}_n \rightarrow 0$. By Lemma
\ref{aszperf} the commutator $R^{\prime}_k$ satisfies
\begin{equation*}
\log |N_n| \ge \log |\RR_k^{\prime}| \ge \log |\Gamma_n|\,
(1-\epsilon^{\prime}_n - \delta_n)
\end{equation*}
where $\delta_n=C \min \{n,-\log_{p} \epsilon^{\prime}\}^{-1/2}$. Letting $%
n\rightarrow \infty$ we get $\gamma_n(N) \rightarrow 1$, as
required.
\end{proof}

We are ready to prove Conjecture \ref{propsub} in the $1$-dimensional
spherically transitive case.

\newtheorem*{traprop.thm}{Theorem \ref{traprop}}
\begin{traprop.thm}
\trapropthm
\end{traprop.thm}

In particular, $G$ involves every finite $p$-group. Recall that
a group $G$ {\bf involves} a group $H$ if it can be obtained as
a quotient of a finite
index subgroup of $G$. Our first lemma confirms this for every $1 $%
-dimensional subgroup.

\begin{lemma}
\lb{involve}Let $G\subseteq \Gamma (p)$ be a closed subgroup of
dimension $1$. Then $G$ involves every finite $p$-group.
\end{lemma}

\begin{proof}
It suffices to prove that $G$ involves $\Gamma _{n}$ for every
$n$. We define the {\bf $n$-sample} $S_{v}$ of $G$ at the
vertex $v$ of depth $k$ as follows. Take the stabilizer of
level $k$. Project this group onto the descendant subtree of
$v$ and cut it at level $k+n$. Let $S_{v}\subseteq \Gamma _{n}$
denote the action on this subtree.

Now for any $\ell$ it is easy to see that $G_{n\ell}=G\Gamma^{(n\ell)}/%
\Gamma^{(n\ell)}$ satisfies
\begin{equation*}
\left| G_{n\ell}\right| \leq \prod_{v}\left| S_v\right|
\end{equation*}
where the product is taken over all vertices in $T_{n\ell}$ at
levels divisible by $n$. Then
\begin{equation*}
\gamma _{n\ell}(G)\leq \frac{\mbox{$\sum_v$}\log \left| S_{v}\right|}{%
\mbox{$\sum_v$}\log \left| \Gamma_n\right|} \leq \max_v {\frac{\log |S_v| }{%
\log |\Gamma_n|}}.
\end{equation*}
As $\ell\rightarrow \infty$, the left hand side converges to
$1$, so there exists a $v$ such that $S_{v}=\Gamma _{n}$. It is
easy to see that $G$ involves all its samples.
\end{proof}

As a corollary, $1$-dimensional groups cannot satisfy any
pro-$p$ identity. It seems plausible that the lemma holds for
arbitrary positive-dimensional closed subgroups.

Let $G\subseteq \Gamma (p)$ be a closed subgroup and let $r$ be an
infinite ray with vertex $r_{n}$ at level $n$. There are two
possible notions of the $n$th stabilizer of $r$ in $G$: the
stabilizer subgroup $F_{n}\subseteq
G_n=G\Gamma^{(n)}/\Gamma^{(n)}$ of $r_n$, and the congruence quotient $%
K_n=K\Gamma^{(n)}/\Gamma^{(n)}$ of the stabilizer $K\subseteq
G$ of $r$. In general, $K_{n}\subseteq F_{n}$ but they need not
coincide. However, equality holds for spherically transitive
groups of high enough dimension.

\begin{lemma}
\lb{tulel}Let $G\subseteq \Gamma (p)$ be a spherically
transitive closed subgroup with $\dimh G>1/p$ and let $r$ be an
infinite ray. Then there exists
$n_{0}$, such that for all $n>n_{0}$ we have $K_{n}=F_{n}$. In particular, $%
\dimh K=\dimh G$.
\end{lemma}

\begin{proof}
Let $V_{n}=(G\cap \Gamma_n)\Gamma_{n+1}/\Gamma_{n+1}$, the
level $n$ subspace of $G$. Then we have $\dim V_{n}\leq p^{n}$
and $\log_p |G_n| = \dim V_0 +\ldots +\dim V_{n-1}$. Our
assumption and a straightforward computation implies that there
exists $n_{0}$ such that for all $n>n_{0}$ the subspace $V_{n}$
is nontrivial.

Now let $n>n_{0}$ and let $g_n\in F_{n}$. We have to show that
$g_n$ can be extended to $K$, that is, there exists $g\in K$ so
that $g / \Gamma^{(n)} =g_n$. In fact, using induction, it
suffices to show that there exists an extension $g_{n+1}$ to
$F_{n+1}$, the stabilizer of the vertex $v$ of $r$ at level
$n+1$.

Let $S$ be the set of possible extensions of $g$ to $G_{n+1}$
Then $S$ forms a coset of $V_{n}$ in $G_{n+1}$. Now $G$ is
spherically transitive and $V_{n} $ is a nontrivial
$G_{n}$-invariant subspace, thus the projection of $V_{n}$ to
the coordinate $v\in r$ is the entire $\mathbb{F}_{p}$. Since
$S$ stabilizes $v$, the projection of $S$ to $v$ also equals
$\mathbb{F}_{p}$. So there exist $g_{n+1}\in S$ which
stabilizes the children of $v$, as required.

For $n>n_{0}$ we have $\left| G_{n}:K_{n}\right| =\left|
G_{n}:F_{n}\right| =p^{n}$, so $\dimh K=\dimh G$.
\end{proof}


\begin{proofof}{Theorem \ref{traprop}}
Let $r$ be the leftmost infinite ray, that is, the set of
vertices indexed by sequences of zeros. Let $F_{k}$ denote the
stabilizer of $r$ in the finite group $\Gamma _{k}(p)$. Let $K$
be the stabilizer of $r$ in $G$. We claim that $K$ has a closed
normal subgroup such that the quotient group is isomorphic to
the full direct product $D=\prod_{k}F_{k}$.

Choose $n_{0}$ according to Lemma \ref{tulel}. From Lemma
\ref{involve} it follows that for every $k$ there exists
infinitely many vertex $v$ such that the $k$-sample
$S_{v}=\Gamma _{k}$. Since $G$ is spherically transitive,
samples at the same level are isomorphic, so we can choose a
strictly increasing sequence $n_{k}$ such that
$S_{v_{k}}=\Gamma _{k}$ where the
vertex $v_{k}$ is at level $n_{k}$ on the ray $r$. We can also assume that $%
n_{k+1}-n_{k}>k$ for all $k$. Now $K$ naturally projects into $\prod_{k}F_{k}
$ by the following function $\pi $. For $g\in K$ let the $k$-th coordinate $%
\pi _{k}(g)$ be the action of $g$ on the subtree of length $k$ starting at
the vertex $v_{k}\in r$. As $g$ stabilizes $r$, the image $\pi _{k}(g)$ will
lie in $F_{k}$.

We show that the projection $\pi $ is surjective, that is, for every
$h\in
\prod_{k}F_{k}$ we can find $g\in K$ such that $\pi (g)=h$. We construct $g$
by a series of level extensions. From level $n_{k}$ to level $n_{k}+k$ we
use the following argument. Since $g$ stabilizes $v_{k}$, the set of
possible extensions of $g$ under $v_{k}$ is the union of cosets of the
sample $S_{v_{k}}(G)$. But $n_{k}$ was chosen in a way such that the $k$%
-sample $S_{v_{k}}(G)$ is the full $\Gamma _{k}(p)$. This implies that up to
level $n_{k}+k$ we can prescribe $g$ in an arbitrary way, so we can assume $%
\pi _{k}(g)=h_{k}$. For levels between $n_{k}+k$ and $n_{k+1}$
we use Lemma \ref{tulel} to obtain an extension which
stabilizes $v_{k+1}$. Since $G$ is closed, our series of
extensions produces an element of $G$, so the claim holds and
$K$ projects onto $D$.

Since $\Gamma _{k-1}(p)$ is a subgroup of $F_{k}$ for all $k$, every finite $%
p$-group can be embedded into $F_{k}$ for large enough $k$.
This implies that the full product $D$ has a nonabelian free
pro-$p$ subgroup topologically generated, say, by the elements $x,\,y$. Then
any two preimages $\tilde{x},\,\tilde{y}$ in $G$ generate a
nonabelian free pro-$p$ subgroup, since pro-$p$ words satisfied
by $\tilde{x}$ and $\tilde{y}$ are also satisfied by $x$ and
$y$. The proof is complete.
\end{proofof}

\cite{MR2001i:20060} showed that just infinite pro-$p$ branch
groups contain dense free subgroups (on the existence of dense
free subgroups in profinite groups see \cite{MR2002g:20114} and
\cite{MR2001h:22017}). The same holds for $1$-dimensional
subgroups. Let $d(G)$ denote the minimal number of topological
generators for $G$ (this might be infinite). An abstract
subgroup of a topological group $G$ is {\bf dense} if its
closure equals $G$.

\newtheorem*{densefree.thm}{Theorem \ref{densefree}}
\begin{densefree.thm}
\densefreethm
\end{densefree.thm}

In other words, the free group $F_{k}$ is residually $S$ where
$S$ denotes the set of congruence quotients of $G$.

\begin{proof}
By Corollary \ref{1dfree} the subgroup generated by $k$ random
elements is free of rank $k$ with probability 1.

If $k=\infty$, then even the set of $k$ random elements is
dense in $G$ with probability $1$, since $G\subseteq \Gamma(p)$
has a countable base. If $d(G)\le k<\infty$ then topological
generation depends only on the Frattini quotient $C_{p}^{d(G)}$
of $G$. This is generated by $k\geq d(G)$ random elements with
positive probability.
\end{proof}

{\noindent \bf Acknowledgments.} We thank Yuval Peres for
recommending the relevant branching random walk references,
P\'{e}ter P\'{a}l P\'{a}lfy for bringing Tur\'{a}n's problem to
our attention, and L\'aszl\'o Pyber for asking a question
leading to Corollary \ref{sokfelold}. We thank Russell Lyons
and P\'{e}ter P\'{a}l P\'{a}lfy for their remarks about
previous versions of this paper.

\bibliography{tg1}

\end{document}